\documentclass[review,final]{siamart1116}

\usepackage{tikz}
\usepackage{amssymb}
\usepackage{amsfonts}
\usepackage{amsmath}
\usepackage{color}
\usepackage{graphicx}
\usepackage{epsfig,mathrsfs}
\usepackage[bf,SL,BF]{subfigure}
\usepackage{fancyhdr}
\usepackage{caption}
\usepackage{wrapfig}
\usepackage{cases}

\DeclareSymbolFont{largesymbol}{OMX}{yhex}{m}{n}
\DeclareMathAccent{\Widehat}{\mathord}{largesymbol}{"62}

\begin{document}
 \pagenumbering{arabic}

\title{{\bf \noindent Boundary problems for the fractional and tempered fractional operators}\thanks{This work was partially supported by NSFC 11421101, 11421110001, 11626250 and 11671182. WD thanks Mark M. Meerschaert and Zhen-Qing Chen for the discussions.}}
\author{
  Weihua Deng\thanks{School of Mathematics and Statistics, Gansu Key Laboratory of Applied Mathematics and Complex Systems, Lanzhou University, Lanzhou 730000, P.R. China. Email: dengwh@lzu.edu.cn}
  \and
  Buyang Li\thanks{Department of Applied Mathematics,
	The Hong Kong Polytechnic University, Hung Hom, Hong Kong. Email: buyang.li@polyu.edu.hk}
  \and
  Wenyi Tian\thanks{Center for Applied Mathematics, Tianjin University, Tianjin 300072, P.R. China. Email: wenyi.tian@tju.edu.cn}
  \and
  Pingwen Zhang\thanks{School of Mathematical Sciences, Laboratory of Mathematics and Applied Mathematics, Peking University, Beijing 100871,  P.R. China. Email: pzhang@pku.edu.cn}
}

\maketitle

\begin{abstract}	
For characterizing the Brownian motion in a bounded domain: $\Omega$, it is well-known that the boundary conditions of the classical diffusion equation just rely on the given information of the solution along the boundary of a domain; on the contrary, for the L\'evy flights or tempered L\'evy flights in a bounded domain, it involves the information of a solution in the complementary set of $\Omega$, i.e., $\mathbb{R}^n\backslash \Omega$, with the potential reason that paths of the corresponding  stochastic process are discontinuous. Guided by probability intuitions and the stochastic perspectives of anomalous diffusion, we show the reasonable ways, ensuring the clear physical meaning and well-posedness of the partial differential equations (PDEs), of specifying `boundary' conditions for space fractional PDEs modeling the anomalous diffusion. Some properties of the operators are discussed, and the well-posednesses of the PDEs with generalized boundary conditions are proved.

\end{abstract}


\section{Introduction}
The phrase `anomalous is normal' says that anomalous diffusion phenomena are ubiquitous in the natural world. It was first used in the title of \cite{SanchoLLSR:2004}, which reveals that the diffusion of classical particles on a solid surface has rich anomalous behaviour controlled by the friction coefficient. In fact,  anomalous diffusion is no longer a young topic. In the review paper \cite{BouchaudG:1990}, the evolution of particles in disordered environments was investigated; the specific effects of a bias on anomalous diffusion were considered; and the generalizations of Einstein's relation in the presence of disorder were discussed. With the rapid development of the study of anomalous dynamics in diverse field, some deterministic equations are derived, governing the macroscopic behaviour of anomalous diffusion. In 2000, Metzler and Klafter published the survey paper \cite{MetzlerK:2000} for the equations governing transport dynamics in complex system with anomalous diffusion and non-exponential relaxation patterns, i.e., fractional kinetic equations of the diffusion,  advection-diffusion, and Fokker-Planck type, derived asymptotically from basic random walk models and a generalized master equation. Many mathematicians have been involved in the research of fractional partial differential equations (PDEs). For fractional PDEs in a bounded domain $\Omega$,
an important question is how to introduce physically meaningful and mathematically well-posed boundary conditions on $\partial \Omega$ or $\mathbb{R}^n\backslash\Omega$.

Microscopically, diffusion is the net movement of particles from a region of higher concentration to a region of lower concentration; for the normal diffusion (Brownian motion), the second moment of the particle trajectories is a linear function of the time $t$; naturally, if it is a nonlinear function of $t$, we call the corresponding diffusion process anomalous diffusion or non-Brownian diffusion  \cite{MetzlerK:2000}. The microscopic (stochastic) models describing anomalous diffusion include continuous time random walks (CTRWs), Langevin type equation, L\'evy processes, subordinated L\'evy processes, and fractional Brownian motion,  etc.. The CTRWs contain two important random variables describing the motion of particles \cite{MontrollW:1965}, i.e., the  waiting time ${\xi}$ and jump length $\eta$. If both the first moment of ${\xi}$ and the second moment of $\eta$ are finite in the scaling limit, then the CTRWs approximate Brownian motion. On the contrary, if one of them is divergent, then the CTRWs characterize anomalous diffusion. Two of the most important CTRW models are L\'evy flights and L\'evy walks. For L\'evy flights, the ${\xi}$ with finite first moment  and $\eta$ with infinite second moment are independent, leading to infinite propagation speed and the divergent second moments of the distribution of the particles. This causes much difficulty in relating the models to experimental data, especially when analyzing the scaling of the measured moments in time \cite{ZaburdaevDK:2015}.
With coupled distribution of ${\xi}$ and $\eta$ (the infinite speed is penalized by the corresponding waiting times), we get the so-called  L\'evy walks \cite{ZaburdaevDK:2015}.  Another idea to ensure that the processes have bounded moments is to truncate the long tailed probability distribution of L\'evy flights \cite{MantegnaS:1994}; they still look like a L\'evy flight in not too long a time.  Currently, the most popular way to do the truncation is to use the exponential tempering, offering the technical advantage of still being an infinitely divisible L\'evy process after the operation \cite{MeerschaertS:2012}. The L\'evy process to describe anomalous diffusion is the scaling limit of CTRWs with independent ${\xi}$ and $\eta$. It is characterized by its characteristic function. Except Brownian motion with drift, the paths of all other proper L\'evy processes are discontinuous. Sometimes, the L\'evy flights are conveniently described by the Brownian motion subordinated to a L\'evy process \cite{ChenS:2005}.
Fractional Brownian motions are often taken as the models to characterize subdiffusion \cite{MandelbrotN:1968}.

Macroscopically, fractional (nonlocal) PDEs are the most popular and effective models for anomalous diffusion, derived from the microscopic models. The solution of fractional PDEs is generally the probability density function (PDF) of the position of the particles undergoing anomalous dynamics; with the deepening of research, the fractional PDEs governing the functional distribution of particles' trajectories are also developed \cite{TurgemanCB:2009, WuDB:2016}. Two ways are usually used to derive the fractional PDEs. One is based on the Montroll-Weiss equation \cite{MontrollW:1965}, i.e., in Fourier-Laplace space, the PDF $p({\bf X},t)$ obeys
\begin{equation} \label{Montroll-WeissEq}
\hat p({\bf k},u)=\frac{1-\phi(u)}{u}\cdot \frac{\hat p_0({\bf k})}{1-\Psi(u,{\bf k})},
\end{equation}
where $\hat p_0({\bf k})$ is the Fourier transform of the initial data; $\phi(u)$ is the Laplace transform of the PDF of waiting times ${\xi}$ and $\Psi(u,{\bf k})$ the Laplace and the Fourier transforms of the joint PDF of waiting times ${\xi}$ and jump length $\eta$. If ${\xi}$ and $\eta$ are independent, then $\Psi(u,{\bf k})=\phi(u)\psi({\bf k})$, where $\psi({\bf k})$ is the Fourier transform of the PDF of $\eta$. Another way is based on the characteristic function of the $\alpha$-stable L\'evy motion, being the scaling limit of the CTRW model with power law distribution of jump length $\eta$. In the high dimensional case, it is more convenient to make the derivation by using the characteristic function of the stochastic process. According to the L\'evy-Khinchin formula \cite{Applebaum:2009}, the characteristic function of L\'evy process has a specific form
\begin{equation} \label{CharacteristicFunction}
\int_{\mathbb{R}^n} e^{i{\bf k}\cdot {\bf X}}p({\bf X},t) {\bf dX}={\bf E} (e^{i {\bf k}\cdot {\bf X} })= e^{t\Phi({\bf k})},
\end{equation}
where
$$
\Phi({\bf k})=i{\bf a \cdot k}-\frac{1}{2} ({\bf k \cdot bk})+\int_{\mathbb{R}^n \backslash \{0\}} \left[e^{i \bf k \cdot X }-1-i({\bf k \cdot X}) \chi_{\{|\bf{X}|<1\}} \right] \nu(d \bf{X});
$$
here $\chi_I$ is the indicator function of the set $I$, ${\bf a} \in \mathbb{R}^n$, ${\bf b}$ is a positive definite symmetric $n \times n$ matrix and $\nu$ is a sigma-finite L\'evy measure on $\mathbb{R}^n \backslash \{0\}$. When ${\bf a}$ and ${\bf b}$ are zero and
\begin{equation}\label{PowerLawMeasure}
\nu(d {\bf X})=\frac{\beta \Gamma(\frac{n+\beta}{2})}{2^{1-\beta} \pi^{n/2}\Gamma(1-\beta/2)}|{\bf X}|^{-\beta-n}{\bf dX},
\end{equation}
the process is a rotationally symmetric $\beta$-stable L\'evy motion and its PDF solves
\begin{equation}\label{FractionalLaplacian}
\frac{\partial p({\bf X},t)}{\partial t}=\Delta^{\beta/2}p({\bf X},t),
\end{equation}
where $\mathcal{F}(\Delta^{\beta/2}p({\bf X},t))=-|{\bf k}|^\beta \mathcal{F}(p({\bf X},t))$ \cite{Silvestre:2007}. If replacing   (\ref{PowerLawMeasure}) by the measure of isotropic tempered power law with the tempering exponent $\lambda$, then we get the corresponding PDF evolution equation
\begin{equation}\label{TemperedFractionalLaplacian}
\frac{\partial p({\bf X},t)}{\partial t}=(\Delta+\lambda)^{\beta/2}p({\bf X},t),
\end{equation}
where $(\Delta+\lambda)^{\beta/2}$ is defined by (\ref{TemperedFracLaplOperator}) in physical space and by (\ref{FourierTemperedLaplacian}) in Fourier space.

%

\begin{figure}
	\begin{center}
		\includegraphics[width=4.2in,height=1.3in,angle=0]{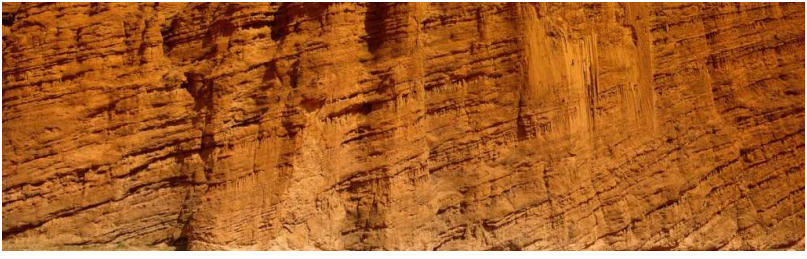}
		\caption{Sketch map for the physical environment suitable for Eq. (\ref{HVequation1}). }\label{Figure1}
	\end{center}
\end{figure}

In practice, the choice of $\nu(d {\bf X})$ depends strongly on the concrete physical environment. For example, Figure \ref{Figure1} clearly shows the horizontal and vertical structure. So, we need to take the measure as (if it is superdiffusion)
\begin{equation} \label{HVmeasure}
\begin{aligned}
\nu(d {\bf X})&=\nu(d {\bf x}_1 d{\bf x}_2) =
\frac{\beta_1 \Gamma(\frac{1+\beta_1}{2})}{2^{1-\beta_1} \pi^{1/2}\Gamma(1-\beta_1/2)}|{\bf x}_1|^{-\beta_1-1} \delta({\bf x}_2){ d{\bf x}_1 d{\bf x}_2}
\\
& ~~ ~~~ ~~~ ~~~ ~~~ ~~~ ~~~  +
\frac{\beta_2 \Gamma(\frac{1+\beta_2}{2})}{2^{1-\beta_2} \pi^{1/2}\Gamma(1-\beta_2/2)}\delta({\bf x}_1)|{\bf x}_2|^{-\beta_2-1}{d{\bf x}_1 d{\bf x}_2},
\end{aligned}
\end{equation}
where $\beta_1$ and $\beta_2$ belong to $(0,2)$. If ${\bf a}$ and ${\bf b}$ equal to zero, then it leads to diffusion equation
\begin{equation} \label{HVequation1}
 \frac{\partial p({\bf x}_1,{\bf x}_2,t)}{\partial t}=\frac{\partial^{\beta_1} p({\bf x}_1,{\bf x}_2,t)}{\partial |{\bf x}_1|^{\beta_1} }+\frac{\partial^{\beta_2}  p({\bf x}_1,{\bf x}_2,t)}{\partial |{\bf x}_2|^{\beta_2}}.
\end{equation}

Under the guidelines of probability intuitions and stochastic perspectives \cite{GuanM:2005} of L\'evy flights or tempered L\'evy flights, we discuss the reasonable ways of defining fractional partial differential operators and specifying the `boundary' conditions for their macroscopic descriptions, i.e., the PDEs of the types Eqs. (\ref{FractionalLaplacian}), (\ref{TemperedFractionalLaplacian}), (\ref{HVequation1}), and their extensions, e.g., the fractional Feynman-Kac equations \cite{TurgemanCB:2009,WuDB:2016}. For the related discussions on the nonlocal diffusion problems from a mathematical point of view, one can see the review paper \cite{DuGLZ:2012}. The divergence of the second moment and the discontinuity of the paths of L\'evy flights predicate that the corresponding diffusion operators should defined on $\mathbb{R}^n$, which further signify that if we are solving the equations in a bounded domain $\Omega$, the information in $\mathbb{R}^n \backslash \Omega$ should also be involved. We will show that the generalized Dirichlet type boundary conditions should be specified as $p({\bf X},t)|_{\mathbb{R}^n\backslash \Omega}=g({\bf X},t)$. If the particles are killed after leaving the domain $\Omega$, then $g({\bf X},t) \equiv 0$, i.e., the so-called absorbing boundary conditions. Because of the discontinuity of the jumps of L\'evy flights, a particular concept `escape probability' can be introduced, which means the probability that the particle jumps from the domain $\Omega$ into a domain $H \subset \mathbb{R}^n\backslash\Omega$; for solving the escape probability, one just needs to specify $g({\bf X})=1$ for ${\bf X} \in H$ and $0$ for ${\bf X} \in ({\mathbb{R}}^n\backslash \Omega)\backslash H$ for the corresponding time-independent PDEs. As for the generalized Neumann type boundary conditions, our ideas come from the fact that the continuity equation (conservation law) holds for any kinds of diffusion, since the particles can not be created or destroyed. Based on the continuity equation and the governing equation of the PDF of L\'evy or tempered L\'evy flights, the corresponding flux ${\bf j}$ can be obtained. So the generalized reflecting boundary conditions should be ${\bf j}|_{\mathbb{R}^n \backslash \Omega} \equiv 0$, which implies $(\nabla \cdot {\bf j})|_{\mathbb{R}^n \backslash \Omega} \equiv 0$. Then, the generalized Neumann type boundary conditions are given as $(\nabla \cdot {\bf j})|_{\mathbb{R}^n \backslash \Omega}=g({\bf X},t)$, e.g., for (\ref{FractionalLaplacian}), it should be taken as  $\left(\Delta^{\beta/2}p({\bf X},t)\right)|_{\mathbb{R}^n \backslash \Omega}=g({\bf X},t)$.
The well-posednesses of the equations under our specified generalized Dirichlet or Neumann type boundary   conditions are well established.

Overall, this paper focuses on introducing physically reasonable boundary constraints for a large class of fractional PDEs,
building a bridge between the physical and mathematical communities for studying anomalous diffusion and fractional PDEs.
In the next section, we recall the derivation of fractional PDEs. Some new concepts are introduced, such as the tempered fractional Laplacian, and some properties of anomalous diffusion are found.
In Sec. 3, we discuss the reasonable ways of specifying the generalized boundary conditions for the fractional PDEs governing the position or functional distributions of L\'evy flights and tempered L\'evy flights. In Sec. 4, we prove well-posedness of the fractional PDEs under the generalized Dirichlet and Neumann boundary conditions defined on the complement of the bounded domain.
Conclusion and remarks are given in the last section.

\section{Preliminaries}
For well understanding and inspiring the ways of specifying the `boundary constrains' to PDEs governing the PDF of L\'evy flights or tempered L\'evy flights, we will show the ideas of deriving the microscopic and macroscopic models.

\subsection{Microscopic models for anomalous diffusion}
\begin{figure}
	\includegraphics[scale=0.30]{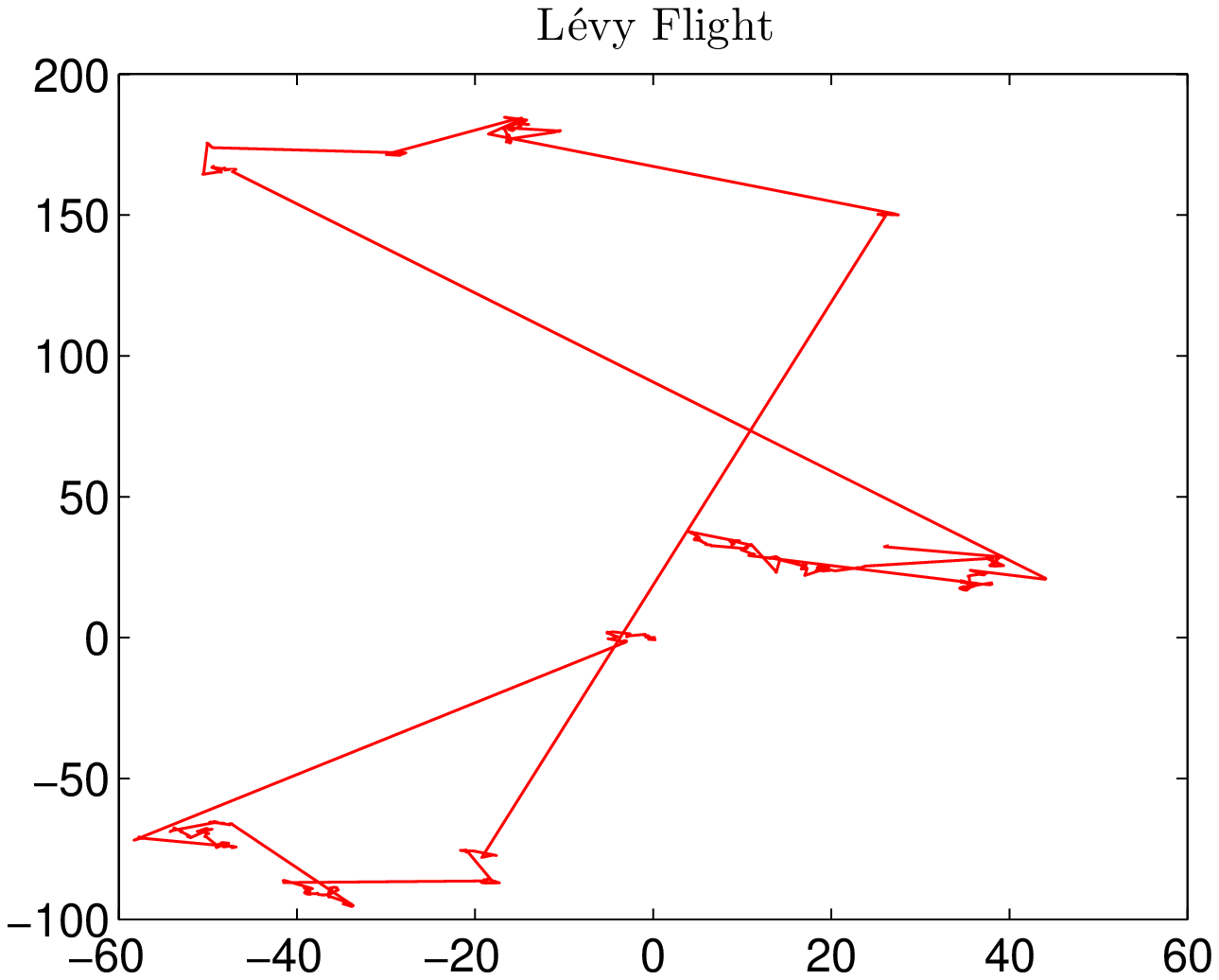}%
	\includegraphics[scale=0.30]{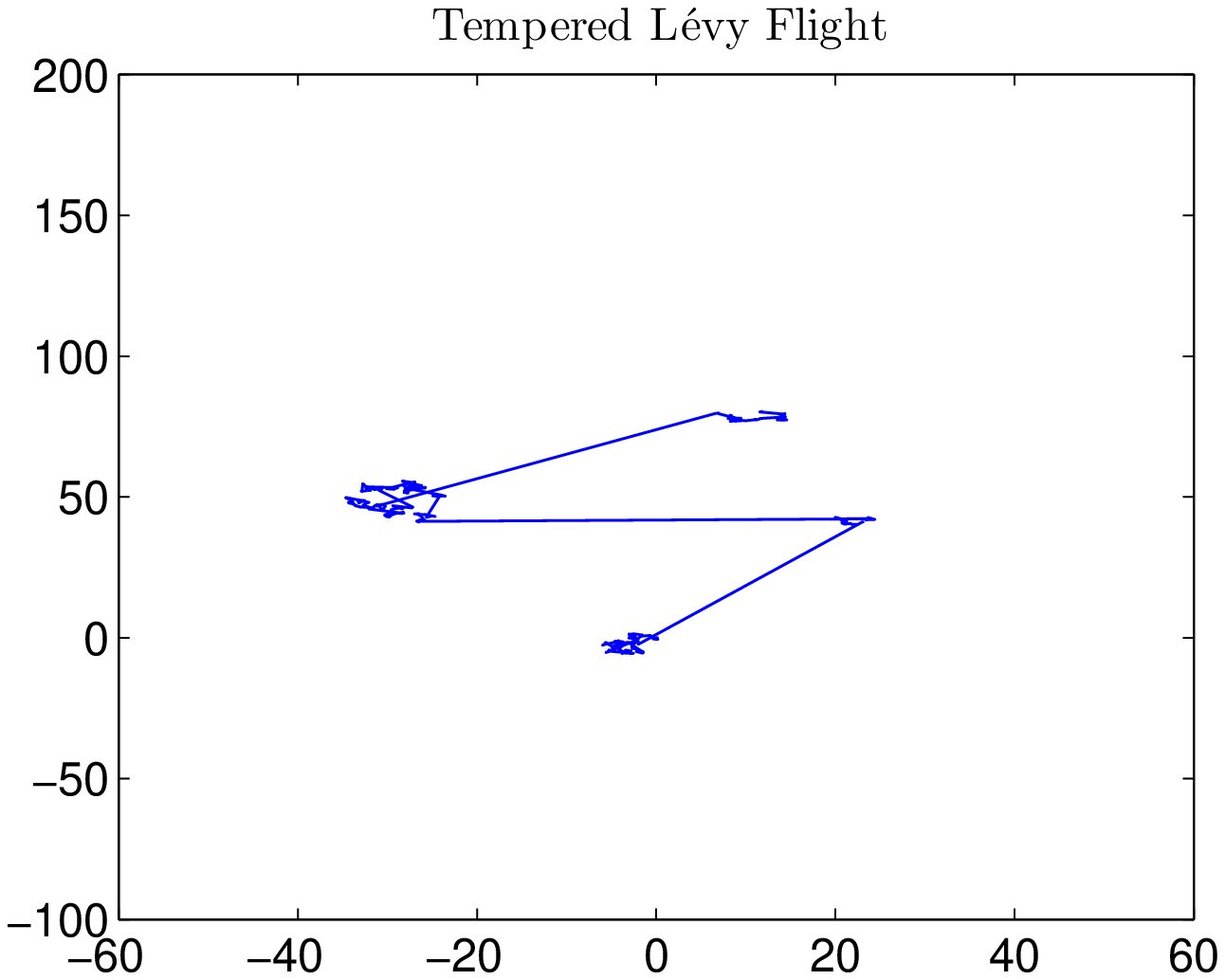}%
	\includegraphics[scale=0.30]{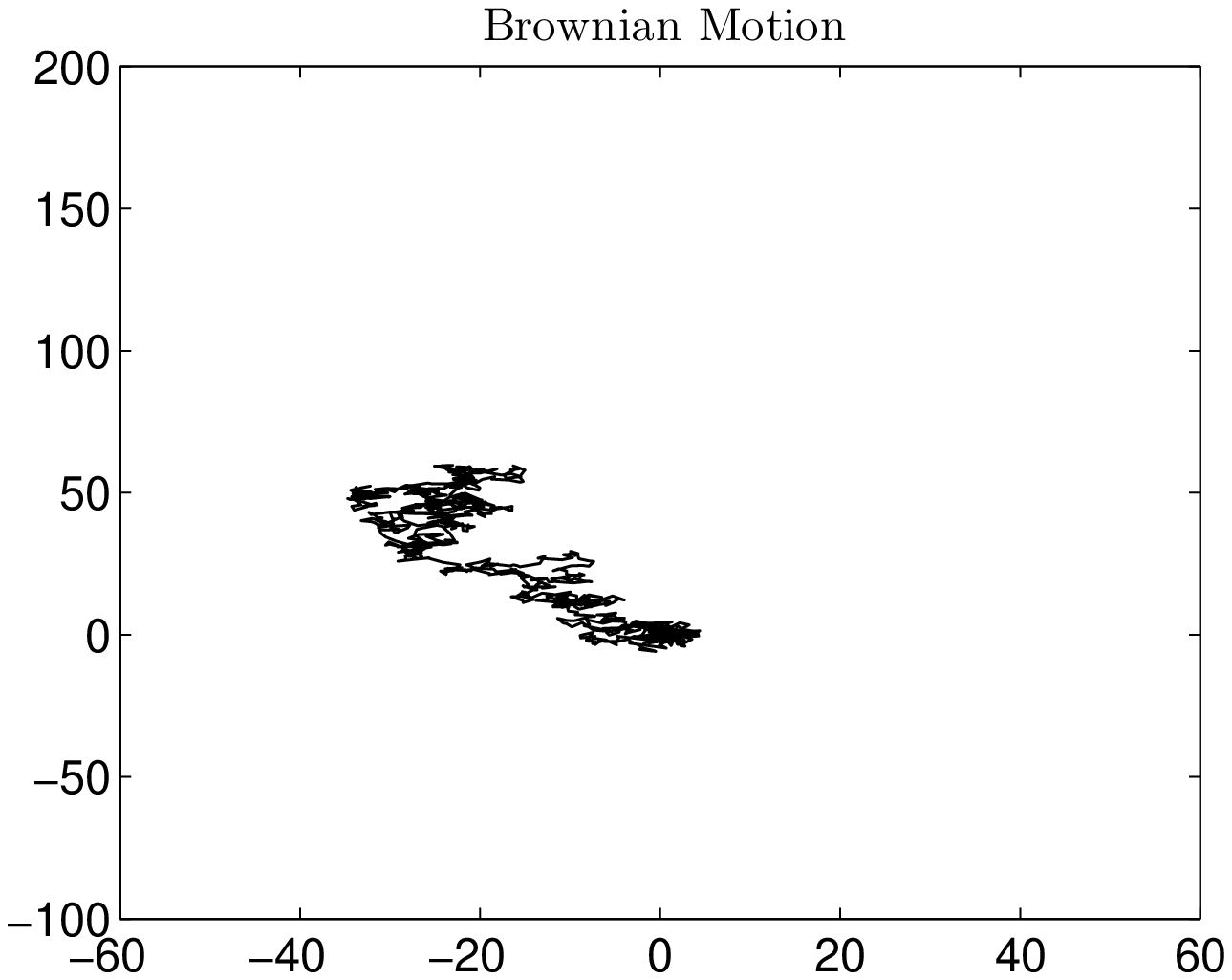}
	\caption{Random trajectories (1000 steps) of L\'evy flight ($\beta=0.8$), tempered L\'evy flight ($\beta=0.8,\,\lambda=0.2$), and Brownian motion.}  \label{Figure2}
\end{figure}

For the microscopic description of the anomalous diffusion, we consider the trajectory of a particle or a stochastic process, i.e., ${\bf X}(t)$. If $\left< |{\bf X}(t)|^2 \right> \sim t$, the process is normal, otherwise it is abnormal. The anomalous diffusions of most often happening in natural world are the cases that  $\left< |{\bf X}(t)|^2 \right> \sim t^\gamma$ with $\gamma\in [0,1) \cup (1,2]$. A L\'evy flight is a random walk in which the jump length has a heavy tailed (power law) probability distribution, i.e., the PDF of jump length $r$ is like $r^{-\beta-n}$ with $\beta \in (0,2)$, and the distribution in direction is uniform.
With the wide applications  of L\'evy flights in characterizing long-range interactions \cite{BarkaiNVBK:2003} or a nontrivial ``crumpled'' topology of a phase (or configuration) space of polymer systems \cite{SokolovMB:1997}, etc, its second and higher moments are divergent, leading to the difficulty in relating models to experimental data.  In fact, for L\'evy flights $\left< |{\bf X}(t)|^\delta \right> \sim t^{\delta/\beta}$ with $0<\delta<\beta \leq 2$. Under the framework of CTRW, the model L\'evy walk \cite{ShlesingerKW:1986} can circumvent this obstacle by putting a larger time cost to a longer displacement, i.e., using the space-time coupled jump length and waiting time distribution $\Psi(r,t)=\frac{1}{2} \delta(r-vt) \phi(t)$.  Another popular model is the so-called tempered L\'evy flights \cite{Koponen:1995}, in which the extremely long jumps is exponentially cut by using the distribution of jump length $e^{-r\lambda}r^{-\beta-n}$ with $\lambda$ being a small modulation parameter (a smooth exponential regression towards  zero). In not too long a time, the tempered L\'evy flights display the dynamical behaviors of L\'evy flights, ultraslowly converging to the normal diffusion. Figure \ref{Figure2} shows the trajectories of $1000$ steps of L\'evy flights, tempered L\'evy flights, and Brownian motion in two dimensions; note the presence of rare but large jumps compared to the Brownian motion, playing the dominant role in the dynamics.

Using Berry-Ess\'een theorem \cite{Feller:1966}, first established in 1941, which applies to the convergence to a Gaussian for a symmetric random walk whose jump probabilities have a finite third moment, we have that for the one dimensional tempered L\'evy flights with the distribution of jump length $C e^{-r\lambda}r^{-\beta-1}$  the convergence speed is
$$
\frac{5}{2\sqrt{2C}}\frac{\Gamma(3-\beta)}{\Gamma(2-\beta)^{3/2}}\lambda^{-\frac{1}{2}\beta}\frac{1}{\sqrt{m}},
$$
which means that the scaling
law for the number of steps needed for Gaussian behavior
to emerge as
\begin{equation} \label{Crossoverlength}
m \sim \lambda^{-\beta}.
\end{equation}
More concretely, letting ${\bf X}_1$, ${\bf X}_2$, $\cdots$, ${\bf X}_m$ be  i.i.d. random variables with PDF $C e^{-r\lambda}r^{-\beta-1}$ and $E(|{\bf X}_1|^2)=\sigma^2>0$, then the cumulative distribution function (CDF) $Q_m$ of ${\bf Y}_m=({\bf X}_1+{\bf X}_2+\cdots+{\bf X}_m)/(\sigma \sqrt{m})$ converges to the CDF $Q({\bf X})$ of the standard normal distribution as
\begin{equation*}
|Q_m({\bf X})-Q({\bf X})|<\frac{5}{2}\frac{\langle|{\bf X}|^3\rangle}{\langle |{\bf X}|^2\rangle^{3/2}}\frac{1}{\sqrt{m}}
=\frac{5}{2\sqrt{2C}}\frac{\Gamma(3-\beta)}{\Gamma(2-\beta)^{3/2}}\lambda^{-\frac{1}{2}\beta}\frac{1}{\sqrt{m}},
\end{equation*}
since
\begin{equation*}
\langle|{\bf X}|^3\rangle=C\int_{-\infty}^{\infty}|{\bf X}|^3e^{-\lambda |{\bf X}|}|{\bf X}|^{-\beta-1}d|{\bf X}|=2C\int_0^{\infty}e^{-\lambda |{\bf X}|}|{\bf X}|^{3-\beta-1}d|{\bf X}|=2C\lambda^{\beta-3}\Gamma(3-\beta)
\end{equation*}
and
\begin{equation*}
\langle|{\bf X}|^2\rangle=C\int_{-\infty}^{\infty}|{\bf X}|^2e^{-\lambda |{\bf X}|}|{\bf X}|^{-\beta-1}d|{\bf X}|=2C\int_0^{\infty}e^{-\lambda |{\bf X}|}|{\bf X}|^{2-\beta-1}d|{\bf X}|=2C\lambda^{\beta-2}\Gamma(2-\beta).
\end{equation*}
From Eq. (\ref{Crossoverlength}), it can be seen that with the decrease of $\lambda$, the required $m$ for the crossover between L\'evy flight behavior and Gaussian behavior increase rapidly. A little bit counterintuitive observation is that the number of variables required to the crossover increases with the increase of $\beta$.

We have described the distributions of jump length for L\'evy flights and tempered L\'evy flights, in which Poisson process is taken as the renewal process. We denote the Poisson process with rate $\zeta>0$ as $N(t)$ and its waiting time distribution between two events is $\zeta e^{-\zeta t}$. Then the  L\'evy flights or tempered L\'evy flights are the compound Poisson process defined as
$
{\bf X}(t)=\sum\limits_{j=0}^{N(t)} {\bf X}_j,
$
where ${\bf X}_j$ are i.i.d. random variables with the distribution of power law or tempered power law. The characteristic function of ${\bf X}(t)$ can be calculated as follows. For real ${\bf k}$, we have
\begin{equation} \label{CharactFuncLevyFlight}
\begin{aligned}
\hat p({\bf k},t) &= {\bf E} (e^{i {\bf k} \cdot {\bf X}(t)}) \\
        &= \sum\limits_{j=0}^\infty {\bf E} (e^{i {\bf k} \cdot {\bf X}(t)}\,|\, N(t)=j)P(N(t)=j) \\
        &= \sum\limits_{j=0}^\infty {\bf E} (e^{i {\bf k} \cdot ({\bf X}_0+{\bf X}_1+\cdots+{\bf X}_j)}\,|\, N(t)=j)P(N(t)=j) \\
        &= \sum\limits_{j=0}^\infty \Phi_0({\bf k})^j \frac{(\zeta t)^j}{j!}e^{-\zeta t} \\
        &= e^{\zeta t(\Phi_0({\bf k})-1)},
\end{aligned}
\end{equation}
where $\Phi_0({\bf k})={\bf E} (e^{i {\bf k}\cdot {\bf X}_0})$, being also the characteristic function of ${\bf X_1}$, ${\bf X_2}$, $\cdots$, ${\bf X_j}$ since they are i.i.d.

In the CTRW model describing one dimensional L\'evy flights or tempered L\'evy flights, the PDF of waiting times is taken as $\zeta e^{-\zeta t}$ with its Laplace transform $\zeta/(u+\zeta)$  and the PDF of jumping length is $c^{-\beta} r^{-\beta-1}$ or $e^{-\lambda r} r^{-\beta-1}$ with its Fourier transform $1-c^\beta |k|^\beta$ or $1-c_{\beta,\lambda} [(\lambda+ik)^\beta-\lambda^\beta]-c_{\beta,\lambda} [(\lambda-ik)^\beta-\lambda^\beta]
$.
Substituting them into the Montroll-Weiss Eq. (\ref{Montroll-WeissEq}) with $\hat p_0(k)=1$ (the initial position of particles is at zero), we get that $\hat p(k,u)$ of L\'evy flights solves
\begin{equation} \label{FLPku}
\hat p(k,u)=\frac{1}{u+\zeta c^\beta |k|^\beta};
\end{equation}
and the $\hat p(k,u)$ of tempered L\'evy flights obeys
\begin{equation} \label{FLTPku}
\hat p(k,u)=\frac{1}{u+\zeta C_{\beta,\lambda} [(\lambda+ik)^\beta-\lambda^\beta]+\zeta C_{\beta,\lambda} [(\lambda-ik)^\beta-\lambda^\beta]}.
\end{equation}
If the subdiffusion is involved, we need to choose the PDF of waiting times as $\tilde{c}^{1+\alpha}t^{-\alpha-1}$ with $\alpha \in (0,1)$ and its Laplace transform $1-\tilde{c}^\alpha u^\alpha$. Then from (\ref{Montroll-WeissEq}), we get that
\begin{equation} \label{FLTPkuT}
\hat p(k,u)=\frac{\tilde{c}^\alpha}{u^{1-\alpha}(1-(1-\tilde{c}_\alpha u^\alpha)\psi(k))}.
\end{equation}

For high dimensional case, the L\'evy flights can also be characterized by Brownian motion subordinated to a L\'evy process. Let ${\bf Y}(t)$ be a Brownian motion with Fourier exponent $-|{\bf k}|^2$ and $S(t)$ a subordinator with Laplace exponent $u^{\beta/2}$ that is independent of ${\bf Y}(t)$. The process ${\bf X}(t)={\bf Y}({S(t)})$ is describing L\'evy flights with Fourier exponent $-|{\bf k}|^\beta$, being the subordinate process of ${\bf Y}(t)$. In effect, denote the characteristic function of ${\bf Y}(t)$ as $\Phi_y({\bf k})$ and the one of $S(t)$ as $\Phi_s(u)$. Then the characteristic function of ${\bf X}(t)$ is as follows:
\begin{equation}\label{SubordinatedProcess}
\begin{aligned}
\hat p_x({\bf k},t) &= \int_{\mathbb{R}^n} e^{i \bf{k} \cdot {\bf X}}p_x( {\bf X},t) d {\bf X} \\
&= \int_0^\infty \int_{\mathbb{R}^n}   e^{i \bf{k} \cdot {\bf Y}}p_y({\bf Y},\tau)d {\bf Y} ~ p_s(\tau,t) d \tau
\\
&=
\int_0^\infty e^{-\tau (- \Phi_y(\bf k))} p_s(\tau,t)d\tau
\\
&=
e^{-t \Phi_s (-\Phi_y(\bf k))},
\end{aligned}
\end{equation}
where $p_x$, $p_y$, and $p_s$, are respectively the PDFs of the stochastic processes ${\bf X}$, ${\bf Y}$, and $S$. Similarly, in the following, we denote $p$ with subscript (lowercase letter) as the PDF of the corresponding stochastic process (uppercase letter).

This paper mainly focuses on L\'evy flights and tempered L\'evy flights. If one is interested in subdiffusion, instead of Poisson process, the fractional Poisson process should be taken as the renewal process, in which the time interval between each pair of events follows the power law distribution. Let ${\bf Y}(t)$ be a general L\'evy process with Fourier exponent $\Phi_y({\bf k})$ and $S(t)$ a strictly increasing subordinator with Laplace exponent $u^\alpha$ ($\alpha \in (0,1)$). Define the inverse subordinator $E(t)=\inf\{ \tau>0:\, S(\tau)>t \}$. Since $t=S(\tau)$ and $\tau=E(t)$ are inverse processes, we have $P(E(t) \le \tau)=P(S(\tau) \ge t)$. Hence
\begin{equation} 
p_e(\tau,t)=\frac{\partial P(E(t) \le \tau)}{\partial \tau}=\frac{\partial}{\partial \tau} \left[1-P(S(\tau) < t) \right]=-\frac{\partial }{\partial \tau} \int_0^t p_s(y,\tau)dy.
\end{equation}
In the above equation, taking Laplace transform w.r.t $t$ leads to
\begin{equation} 
p_e(\tau,u)=-\frac{\partial}{\partial \tau} u^{-1}e^{-\tau u^\alpha}=u^{\alpha-1}e^{-\tau u^\alpha}.
\end{equation}
For the PDF $p_x({\bf X},t)$ of  ${\bf X}(t)={\bf Y}(E(t))$, there holds
\begin{equation}\label{PDFRel}
p_x({\bf X},t)=\int_0^\infty p_y({\bf X},\tau) p_e(\tau,t) d \tau.
\end{equation}
Performing Fourier transform w.r.t. ${\bf X}$ and Laplace transform w.r.t. $t$ to the above equation results in
\begin{equation}\label{PDFRel2}
\begin{aligned}
\hat p_x({\bf k}, u) &= \int_0^\infty \hat p_y({\bf k},\tau)p_e(\tau,u)d \tau
\\
&= \int_0^\infty e^{-\tau \Phi_y({\bf k})} u^{\alpha-1}e^{-\tau u^\alpha}d \tau
\\
&=
\frac{u^{\alpha-1}}{u^\alpha+\Phi_y({\bf k})}.
\end{aligned}
\end{equation}
{\bf Remark.} According to Fogedby \cite{Fogedby:1994}, the stochastic trajectories of (scale limited) CTRW ${\bf X} (E_t)$ can also be expressed in terms of the coupled Langevin equation
\begin{equation}\label{LangevinEq}
\left\{
\begin{aligned}
&\dot{\bf X}(\tau)=F({\bf X}(\tau))+\eta(\tau), \\
&\dot{S}(\tau)=\xi(\tau),
\end{aligned}
\right.
\end{equation}
where $F({\bf X})$ is a vector field; $E_t$ is the inverse process of $S(t)$; the noises $\eta(\tau)$ and $\xi(\tau)$ are statistically independent, corresponding to the distributions of jump length and waiting times.

\subsection{Derivation of the macroscopic description from the microscopic models}

This section focuses on the derivation of the deterministic equations governing the PDF of position of the particles undergoing anomalous diffusion. It shows that the operators related to (tempered) power law jump lengths should be defined on the whole unbounded domain $\mathbb{R}^n$, which can also be inspired by the rare but extremely long jump lengths displayed in Figure \ref{Figure2}; the fact that among all proper L\'evy processes Brownian motion is the unique one with continuous paths further consolidates the reasonable way of defining the operators.
We derive the PDEs based on Eqs. (\ref{CharactFuncLevyFlight}), (\ref{SubordinatedProcess}), and  (\ref{PDFRel}), since they apply for both one and higher dimensional cases. For one dimensional case, sometimes it is convenient to use (\ref{FLPku}), (\ref{FLTPku}), and  (\ref{FLTPkuT}).

When the diffusion process is rotationally symmetric $\beta$-stable, i.e., it is isotropic with PDF of jump length $c_{\beta,n}r^{-\beta-n}$ and its Fourier transform $1-{|\bf k|}^\beta$, where $n$ is the space dimension.  In Eq. (\ref{CharactFuncLevyFlight}), taking $\zeta$ equal to $1$, we get the Cauchy equation
\begin{equation} \label{CauchyEq}
\frac{d \hat p({\bf k},t)}{dt}=-|{\bf k}|^\beta \hat p({\bf k},t).
\end{equation}

Performing inverse Fourier transform to the above equation leads to
\begin{equation} \label{FracLapl}
\frac{\partial p({\bf X},t)}{\partial t}=\Delta^{\beta/2}p({\bf X},t),
\end{equation}
where
\begin{equation} \label{FracLaplOperator}
\begin{aligned}
\Delta^{\beta/2}p({\bf X},t) &= -c_{n,\beta} \lim\limits_{\varepsilon \rightarrow 0^+} \int_{\mathscr C B_\varepsilon ({\bf X})} \frac{p({\bf X},t)-p({\bf Y},t)}{|{\bf X}- {\bf Y}|^{n+\beta}}d{\bf Y} \\
\\ &=
\frac{1}{2} c_{n,\beta} \int_{\mathbb{R}^n} \frac{p({\bf X}+{\bf Y},t)+p({\bf X}-{\bf Y},t)-2 \cdot p({\bf X},t)}{|{\bf Y}|^{n+\beta}}d{\bf Y}
\end{aligned}
\end{equation}
with \cite{DiNezzaPV:2012}
\begin{equation} \label{CoefficientLaplacian}
c_{n,\beta}=\frac{\beta \Gamma(\frac{n+\beta}{2})}{2^{1-\beta} \pi^{n/2}\Gamma(1-\beta/2)}.
\end{equation}
For the more general cases of Eq. (\ref{CharactFuncLevyFlight}), there is the Cauchy equation
\begin{equation}\label{GenCauchyEq}
\frac{d \hat p({\bf k},t)}{dt}=(\Phi_0({\bf k})-1) \hat p({\bf k},t),
\end{equation}
so the PDF of the stochastic process ${\bf X}$ solves (taking $\zeta=1$)
\begin{equation} \label{PDFeq}
\begin{aligned}
\frac{\partial p({\bf X},t)}{\partial t} &= \mathcal{F}^{-1} \{(\Phi_0({\bf k})-1) \hat p({\bf k},t) \}
\\
&= \int_{{\mathbb R}^n\backslash \{0\}} [p({\bf X}+{\bf Y},t)-p({\bf X},t)] \nu (d {\bf Y}),
\end{aligned}
\end{equation}
where $\nu (d {\bf Y})$ is the probability measure of the jump length. Sometimes, to overcome the possible divergence of the terms on the right hand side of Eq. (\ref{PDFeq}) because of the possible strong singularity of $\nu (d {\bf Y})$ at zero, the term
$$
\Phi_0({\bf k})-1=\int_{{\mathbb R}^n \backslash \{0\}} \left [e^{i {\bf k \cdot Y} }-1  \right] \nu (d {\bf Y})
$$
is approximately replaced by
\begin{equation} \label{TermExtension}
\int_{{\mathbb R}^n \backslash \{0\}}\left [e^{i {\bf k \cdot Y} }-1 - i ({\bf k \cdot Y})_{\chi_{ \{|{\bf Y}|<1\}}} \right] \nu (d {\bf Y});
\end{equation}
then the corresponding modification to Eq. (\ref{PDFeq}) is
\begin{equation} \label{PDFeqExtension}
\frac{\partial p({\bf X},t)}{\partial t} =\int_{{\mathbb R}^n\backslash \{0\}} \left[p({\bf X}+{\bf Y},t)-p({\bf X},t)-\sum\limits_{i=1}^n {\bf y}_i(\partial_i p({\bf X},t))_{\chi_{ \{|{\bf Y}|<1\}}} \right]\nu (d {\bf Y}),
\end{equation}
where ${\bf y}_i$ is the component of ${\bf Y}$, i.e., ${\bf Y}=\{{\bf y}_1, {\bf y}_2, \cdots, {\bf y}_n \}^T$. If $\nu (-d {\bf Y})=\nu (d {\bf Y})$, the integration of the summation term of above equation equals to zero.

If the diffusion is in the environment having a structure like Figure \ref{Figure1}, the probability measure should be taken as
\begin{equation} \label{HVmeasure2}
\begin{aligned}
\nu(d {\bf X}) &= \nu(d {\bf x}_1 d{\bf x}_2 d{\bf x}_3 \cdots d{\bf x}_n)
\\
&=
\frac{\beta_1 \Gamma(\frac{1+\beta_1}{2})}{2^{1-\beta_1} \pi^{1/2}\Gamma(1-\beta_1/2)}|{\bf x}_1|^{-\beta_1-1} \delta({\bf x}_2) \delta({\bf x}_3) \cdots \delta({\bf x}_n) d {\bf x}_1 d{\bf x}_2 d{\bf x}_3 \cdots d{\bf x}_n
\\
&
~~~ +
\frac{\beta_2 \Gamma(\frac{1+\beta_2}{2})}{2^{1-\beta_2} \pi^{1/2}\Gamma(1-\beta_2/2)}|{\bf x}_2|^{-\beta_2-1} \delta({\bf x}_1) \delta({\bf x}_3) \cdots \delta({\bf x}_n) d {\bf x}_1 d{\bf x}_2 d{\bf x}_3 \cdots d{\bf x}_n
+ \cdots\\
&
~~~+
\frac{\beta_n \Gamma(\frac{1+\beta_n}{2})}{2^{1-\beta_n} \pi^{1/2}\Gamma(1-\beta_n/2)}|{\bf x}_n|^{-\beta_n-1} \delta({\bf x}_1) \delta({\bf x}_2) \cdots \delta({\bf x}_{n-1}) d {\bf x}_1 d{\bf x}_2 d{\bf x}_3 \cdots d{\bf x}_n
,
\end{aligned}
\end{equation}
where $\beta_1,\,\beta_2,\cdots,\,\beta_n$ belong to $(0,2)$.
Plugging Eq. (\ref{HVmeasure2}) into Eq. (\ref{PDFeq}) leads to
\begin{equation} \label{HVequation2}
\frac{\partial p({\bf x}_1,\cdots,{\bf x}_n,t)}{\partial t}=\frac{\partial^{\beta_1} p({\bf x}_1,\cdots,{\bf x}_n,t)}{\partial |{\bf x}_1|^{\beta_1} }+\frac{\partial^{\beta_2}  p({\bf x}_1,\cdots,{\bf x}_n,t)}{\partial |{\bf x}_2|^{\beta_2}}+\cdots+\frac{\partial^{\beta_n}  p({\bf x}_1,\cdots,{\bf x}_n,t)}{\partial |{\bf x}_n|^{\beta_n}},
\end{equation}
where
$$
\mathcal{F} \left( \frac{\partial^{\beta_j}  p({\bf x}_1,\cdots,{\bf x}_n,t)}{\partial |{\bf x}_j|^{\beta_j}}  \right)
=-|{\bf k}_j|^{\beta_j} p({\bf x}_1,\cdots, {\bf x}_{j-1}, {\bf k}_j,{\bf x}_{j+1},\cdots, {\bf x}_n,t)
$$
and $\frac{\partial^{\beta_j}  p({\bf x}_1,\cdots,{\bf x}_n,t)}{\partial |{\bf x}_j|^{\beta_j}}$ in physical space is
defined by (\ref{FracLaplOperator}) with $n=1$; in particular, when $\beta_j \in (1,2)$, it can also be written as
\begin{equation} \label{RieszDerivative}\small
\frac{\partial^{\beta_j}  p({\bf x}_1,\cdots,{\bf x}_n,t)}{\partial |{\bf x}_j|^{\beta_j}}=-\frac{1}{2\cos(\beta_j \pi/2) \Gamma(2-\beta_j)} \frac{\partial^2}{\partial {\bf x}_j^2} \int_{-\infty}^{\infty} |{\bf x}_j-{\bf y}|^{1-\beta_j} p({\bf x}_1,\cdots, {\bf y},\cdots, {\bf x}_n,t)d {\bf y}.
\end{equation}
It should be emphasized here that when characterizing diffusion processes related with L\'evy flights the operators should be defined in the whole space. Another issue that also should be stressed is that when $\beta_1=\beta_2=\cdots=\beta_n=1$, Eq. (\ref{HVequation2}) is still describing the phenomena of anomalous diffusion, including the cases that they belong to $(0,1)$; the corresponding `first' order operator is nonlocal, being different from the classical first order operator, but they have the same energy in the sense that
$$
\begin{aligned}
&\mathcal{F} \left( \frac{\partial  p({\bf x}_1,\cdots,{\bf x}_n,t)}{\partial |{\bf x}_j|}  \right)
\overline{\mathcal{F} \left( \frac{\partial  p({\bf x}_1,\cdots,{\bf x}_n,t)}{\partial |{\bf x}_j|}  \right)
}
\\
&=\mathcal{F} \left( \frac{\partial  p({\bf x}_1,\cdots,{\bf x}_n,t)}{\partial {\bf x}_j}  \right)
\overline{\mathcal{F} \left( \frac{\partial  p({\bf x}_1,\cdots,{\bf x}_n,t)}{\partial {\bf x}_j}  \right)
}\\
&=(k_j)^2 \hat p^2({\bf x}_1,\cdots, {\bf x}_{j-1}, {\bf k}_j,{\bf x}_{j+1},\cdots, {\bf x}_n,t);
\end{aligned}
$$
$$
\begin{aligned}
&\mathcal{F} \left( \Delta^{1/2} p({\bf X},t)\right) \overline{\mathcal{F} \left( \Delta^{1/2} p({\bf X},t)\right)} \\
& = \mathcal{F} \left( \nabla p({\bf X},t)  \right) \cdot \overline{\mathcal{F} \left( \nabla p({\bf X},t) \right)}=|{\bf k}|^2 \hat p^2({\bf k},t),
\end{aligned}
$$
even though $\Delta^{1/2}$ and $\nabla$ are completely different operators, where the notation $\overline{v}$ stands for the complex conjugate of $v$.


If the subdiffusion is involved, the derivation of the macroscopic equation should be based on Eq. (\ref{PDFRel2}). For getting the term related to time derivative, the inverse Laplace transform should be performed on $u^\alpha \hat p({\bf k},u)-u^{\alpha-1}$. Since $\hat p({\bf k},t=0)$ is taken as $1$, there exists
\begin{equation} \label{TimeDerivative}
\mathcal{L}^{-1} (u^\alpha \hat p({\bf k},u)-u^{\alpha-1})=\frac{1}{\Gamma(1-\alpha)}
\int_0^t (t-\tau)^{-\alpha} \frac{\partial \hat p({\bf k},\tau)}{\partial \tau} d \tau,
\end{equation}
which is usually denoted as ${_0^CD_t^\alpha} \hat p({\bf k},t)$, the so-called Caputo fractional derivative. So, if both the PDFs of the waiting time and jump lengths of the stochastic process ${\bf X}$ are power law, the corresponding models can be obtained by replacing $\frac{\partial }{\partial t}$ with ${_0^CD_t^\alpha}$ in Eqs. (\ref{FracLapl}), (\ref{PDFeq}),
(\ref{PDFeqExtension}), and (\ref{HVequation2}). Furthermore, if there is an external force $F(\bf X)$ in the considered stochastic process $\bf X$, we need to add an additional term $\nabla \cdot (
F({\bf X})p({\bf X},t))$ on the right hand side of Eqs. (\ref{FracLapl}), (\ref{PDFeq}),
(\ref{PDFeqExtension}), and (\ref{HVequation2}).

Here we turn to another important and interesting topic: tempered L\'evy flights. Practically it is not easy to collect the value of a function in the unbounded area ${\mathbb R}^n\backslash \Omega$. This is one of the achievements of using tempered fractional Laplacian. It is still isotropic but with PDF of jump length
$c_{\beta,n,\lambda}e^{-\lambda r}r^{-\beta-n}$. The PDF of tempered L\'evy flights solves
\begin{equation}\label{temperedLaplacian}
\frac{\partial p({\bf X},t)}{\partial t}=(\Delta+\lambda)^{\beta/2}p({\bf X},t),
\end{equation}
where
\begin{equation} \label{TemperedFracLaplOperator}
\begin{aligned}
(\Delta+\lambda)^{\beta/2}p({\bf X},t) &= -c_{n,\beta,\lambda} \lim\limits_{\varepsilon \rightarrow 0^+} \int_{\mathscr C B_\varepsilon ({\bf X})} \frac{p({\bf X},t)-p({\bf Y},t)}{e^{\lambda |{\bf X}-{\bf Y}|}|{\bf X}- {\bf Y}|^{n+\beta}}d{\bf Y} \\
\\ &=
\frac{1}{2} c_{n,\beta,\lambda} \int_{\mathbb{R}^n} \frac{p({\bf X}+{\bf Y},t)+p({\bf X}-{\bf Y},t)-2 \cdot p({\bf X},t)}{e^{\lambda |{\bf Y}|}|{\bf Y}|^{n+\beta}}d{\bf Y}
\end{aligned}
\end{equation}
with
\begin{equation} \label{CoefficientTemperedLaplacian}
c_{n,\beta,\lambda}=\frac{-\Gamma(\frac{n}{2})}{2\pi^{\frac{n}{2}}\Gamma(-\beta)}.
\end{equation}
The choice of the constant as the one given in (\ref{CoefficientTemperedLaplacian}) leads to
\begin{equation} \label{FourierTemperedLaplacian}
\mathcal{F} \left((\Delta+\lambda)^{\beta/2}p({\bf X},t)\right)=\big(\lambda^{\beta}-(\lambda^2+|{\bf k}|^2)^{\frac{\beta}{2}}+ O(|{\bf k}|^2) \big) \hat p({\bf k},t
) ~~{\rm with}~~ \beta\in (0,1) \cup (1,2).
\end{equation}
However, if $\lambda=0$, one needs to choose the constant as the one given in (\ref{CoefficientLaplacian}) to make sure
 $\mathcal{F} \left(\Delta^{\beta/2}p({\bf X},t)\right)=-|{\bf k}|^{\beta} \hat p({\bf k},t)$. The reason is as follows.

\begin{equation*}
\begin{aligned}
\mathcal{F} \left((\Delta+\lambda)^{\beta/2}p({\bf X},t)\right) &= \frac{1}{2}c_{n,\beta,\lambda}\int_{\mathbb{R}^n}\frac{e^{i{\bf k}\cdot {\bf Y}}+e^{-i{\bf k}\cdot {\bf Y}}-2}{|{\bf Y}|^{n+\beta}}e^{-\lambda|{\bf Y}|}d{\bf Y} \cdot \mathcal{F}(p({\bf X},t))
\\
\\
&=- c_{n,\beta,\lambda}\int_{\mathbb{R}^n}\frac{1-\cos({\bf k}\cdot {\bf Y})}{|{\bf Y}|^{n+\beta}}e^{-\lambda|{\bf Y}|}d{\bf Y}\cdot\mathcal{F}(p({\bf X},t)).
\end{aligned}
\end{equation*}
For $\beta\in(0,1)\cup(1,2)$, then we have
\begin{equation*}\small
\begin{aligned}
&\int_{\mathbb{R}^n}\frac{1-\cos({\bf k}\cdot {\bf Y})}{e^{\lambda|{\bf Y}|}|{\bf Y}|^{n+\beta}}d{\bf Y}=\int_{\mathbb{R}^n}\frac{1-\cos(|{\bf k}|{\bf y}_1)}{e^{\lambda|{\bf Y}|}|{\bf Y}|^{n+\beta}}d{\bf Y}=|{\bf k}|^{\beta}\int_{\mathbb{R}^n}\frac{1-\cos({\bf x}_1)}{|{\bf X}|^{n+\beta}}e^{-\frac{\lambda}{|{\bf k}|}|{\bf X}|}d{\bf X}\\
&=C|{\bf k}|^{\beta}\int_0^\infty \frac{1}{r^{n+\beta}}e^{-\frac{\lambda}{|{\bf k}|}r}r^{n-1}\Big(\int_0^{\pi}\big(1-\cos(r\cos\theta_1)\big)\sin^{n-2}(\theta_1)d\theta_1\Big) dr\\
&=\frac{1}{(-\beta)(-\beta+1)}C|{\bf k}|^{\beta-2}\lambda^2\int_0^{\infty}e^{-\frac{\lambda}{|{\bf k}|}r}r^{-\beta+1}\Big(\int_0^{\pi}\big(1-\cos(r\cos\theta_1)\big)\sin^{n-2}(\theta_1)d\theta_1\Big) dr\\
&~~~~~-\frac{1}{(-\beta)(-\beta+1)}C|{\bf k}|^{\beta-1}\lambda\int_0^{\infty}e^{-\frac{\lambda}{|{\bf k}|}r}r^{-\beta+1}\Big(\int_0^{\pi}\sin(r\cos\theta_1)\big)\sin^{n-2}(\theta_1)\cos(\theta_1)d\theta_1\Big) dr\\
&~~~~~-\frac{1}{-\beta}C|{\bf k}|^{\beta}\int_0^{\infty}e^{-\frac{\lambda}{|{\bf k}|}r}r^{-\beta}\Big(\int_0^{\pi}\sin(r\cos\theta_1)\big)\sin^{n-2}(\theta_1)\cos(\theta_1)d\theta_1\Big) dr\\
&=C\Gamma(-\beta)\frac{\sqrt{\pi}\Gamma(\frac{n-1}{2})}{\Gamma(\frac{n}{2})}\lambda^{\beta}\bigg[1-{_2}F_1\Big(\frac{2-\beta}{2},\frac{3-\beta}{2};\frac{n}{2};-\frac{|{\bf k}|^2}{\lambda^2}\Big)\\
&~~~~~~~~~~~~~~~~~~~~~~~~~~~~~~~~~~~~-\frac{2-\beta}{n}\frac{|{\bf k}|^2}{\lambda^2}{_2}F_1\Big(\frac{3-\beta}{2},2-\frac{\beta}{2};\frac{n}{2}+1;-\frac{|{\bf k}|^2}{\lambda^2}\Big)\\
&~~~~~~~~~~~~~~~~~~~~~~~~~~~~~~~~~~~~-\frac{1-\beta}{n}\frac{|{\bf k}|^2}{\lambda^2}{_2}F_1\Big(\frac{2-\beta}{2},\frac{3-\beta}{2};\frac{n}{2}+1;-\frac{|{\bf k}|^2}{\lambda^2}\Big)\bigg]\\
&=C\Gamma(-\beta)\frac{\sqrt{\pi}\Gamma(\frac{n-1}{2})}{\Gamma(\frac{n}{2})}\bigg[\lambda^{\beta}-\lambda^{\beta}{_2}F_1\Big(-\frac{\beta}{2},\frac{1-\beta}{2};\frac{n}{2};-\frac{|{\bf k}|^2}{\lambda^2}\Big)\bigg]\\
&=C\Gamma(-\beta)\frac{\sqrt{\pi}\Gamma(\frac{n-1}{2})}{\Gamma(\frac{n}{2})}\bigg[\lambda^{\beta}-\lambda^{\beta}\Big(1+\frac{|{\bf k}|^2}{\lambda^2}\Big)^{\frac{\beta}{2}}{_2}F_1\Big(-\frac{\beta}{2},\frac{n+\beta-1}{2};\frac{n}{2};\frac{|{\bf k}|^2}{\lambda^2+|{\bf k}|^2}\Big)\bigg]\\
&=C\Gamma(-\beta)\frac{\sqrt{\pi}\Gamma(\frac{n-1}{2})}{\Gamma(\frac{n}{2})}\bigg[\lambda^{\beta}-(\lambda^2+|{\bf k}|^2)^{\frac{\beta}{2}}{_2}F_1\Big(-\frac{\beta}{2},\frac{n+\beta-1}{2};\frac{n}{2};\frac{|{\bf k}|^2}{\lambda^2+|{\bf k}|^2}\Big)\bigg],
\end{aligned}
\end{equation*}
where ${_2}F_1$ is the Gaussian hypergeometric function and
\begin{equation*}
C=\Big(\int_0^{\pi}\sin^{n-3}(\theta_2)d\theta_2\Big)\cdots\Big(\int_0^{\pi}\sin(\theta_{n-2})d\theta_{n-2}\Big)\Big(\int_0^{2\pi}d\theta_{n-1}\Big)
=\frac{2\pi^{\frac{n-1}{2}}}{\Gamma(\frac{n-1}{2})}.
\end{equation*}
So
\begin{equation*}
 c_{n,\beta,\lambda}=\frac{-\Gamma(\frac{n}{2})}{2\pi^{\frac{n}{2}}\Gamma(-\beta)}.
\end{equation*}
The PDEs for tempered L\'evy flights or tempered L\'evy flights combined with subdiffusion can be similarly derived, as those done in this section for L\'evy flights or L\'evy flights combined with subdiffusion.
 Here, we present the counterpart of Eq. (\ref{HVequation2}),
\begin{equation} \label{TemperedHVequation2}\small
\frac{\partial p({\bf x}_1,\cdots,{\bf x}_n,t)}{\partial t}=\frac{\partial^{\beta_1,\lambda} p({\bf x}_1,\cdots,{\bf x}_n,t)}{\partial |{\bf x}_1|^{\beta_1,\lambda} }+\frac{\partial^{\beta_2,\lambda}  p({\bf x}_1,\cdots,{\bf x}_n,t)}{\partial |{\bf x}_2|^{\beta_2,\lambda}}+\cdots+\frac{\partial^{\beta_n,\lambda}  p({\bf x}_1,\cdots,{\bf x}_n,t)}{\partial |{\bf x}_n|^{\beta_n,\lambda}},
\end{equation}
where the operator $\frac{\partial^{\beta_j,\lambda}  p({\bf x}_1,\cdots,{\bf x}_j,t)}{\partial |{\bf x}_j|^{\beta_j,\lambda}}$ is defined by taking $\beta=\beta_j$ and $n=1$ in Eq. (\ref{TemperedFracLaplOperator}).
Again, even for the tempered L\'evy flights, all the related operators should be defined on the whole space, because of the very rare but still possible unbounded jump lengths.

All the above derived PDEs are governing the PDF of the position of particles. If one wants to dig out more deep informations of the corresponding stochastic processes, analyzing the distribution of the functional defined by $A=\int_0^t U({\bf X(\tau)})d\tau$ is one of the choices, where $U$ is a prespecified function. Denote the PDF of the functional $A$ and position ${\bf X}$ as $G({\bf X},A,t)$ and the counterpart of $A$ in Fourier space as $q$. Then $\hat G({\bf X},q,t)$ solves \cite{TurgemanCB:2009}
\begin{equation} \label{Feynman-Kac}
\frac{\partial \hat G({\bf X},q,t)}{\partial t}=K_{\alpha,\beta} \Delta^{\beta/2} D_t^{1-\alpha} \hat G({\bf X},q,t)+i q U({\bf X}) \hat G({\bf X},q,t)
\end{equation}
for L\'evy flights combined with subdiffusion;
and \cite{WuDB:2016}
\begin{equation} \label{TemperedFeynman-Kac}
\frac{\partial \hat G({\bf X},q,t)}{\partial t}=K_{\alpha,\beta} (\Delta+\lambda)^{\beta/2} D_t^{1-\alpha} \hat G({\bf X},q,t)+i q U({\bf X}) \hat G({\bf X},q,t)
\end{equation}
for tempered L\'evy flights combined with subdiffusion,
where
$$
D_t^{1-\alpha}\hat G({\bf X},q,t)=\frac{1}{\Gamma(\alpha)} \left[ \frac{\partial}{\partial t} -iqU({\bf X})\right] \int_0^t \frac{e^{i(t-\tau)qU({\bf X})}}{(t-\tau)^{1-\alpha}} \hat G({\bf X},q,\tau) d\tau.
$$
If one is only interested in the functional $A$ (not caring position ${\bf X}$), then $\hat G_{{\bf X}_0}(q,t)$ is, respectively, governed by \cite{TurgemanCB:2009}
\begin{equation} \label{BackwardFeynman-Kac}
\frac{\partial \hat  G_{{\bf X}_0}(q,t)}{\partial t}=K_{\alpha,\beta} D_t^{1-\alpha}\Delta^{\beta/2} \hat G_{{\bf X}_0}(q,t)+i q U({\bf X}) \hat G_{{\bf X}_0}(q,t)
\end{equation}
and \cite{WuDB:2016}
\begin{equation} \label{BackwardTemperedFeynman-Kac}
\frac{\partial \hat G_{{\bf X}_0}(q,t)}{\partial t}=K_{\alpha,\beta} D_t^{1-\alpha}(\Delta+\lambda)^{\beta/2} \hat G_{{\bf X}_0}(q,t)+i q U({\bf X}) \hat G_{{\bf X}_0}(q,t)
\end{equation}
for L\'evy flights and tempered L\'evy flights, combined with subdiffusion; the ${\bf X}_0$ in $\hat G_{{\bf X}_0}(q,t)$ means the initial position of particles, being a parameter.

\section{Specifying the generalized boundary conditions for the fractional PDEs}
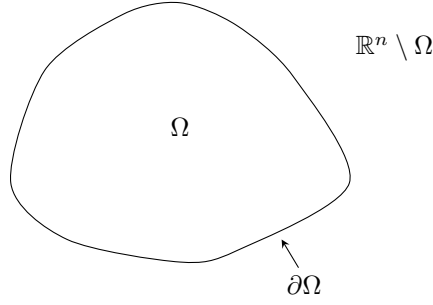
\begin{figure}
	\centering
	\begin{tikzpicture}
	\def\sc{1.5}
	\node at (1.5*\sc,0.5*\sc) {$\Omega$};
	\draw plot[smooth cycle] coordinates %
	{(0,0) (0.3*\sc,1*\sc) (1*\sc,1.5*\sc) (1.5*\sc, 1.6*\sc) (2*\sc,1.4*\sc) (2.5*\sc,0.95*\sc) (3*\sc,0*\sc) (2.0*\sc,-0.6*\sc) (1.5*\sc,-0.7*\sc) (0.5*\sc,-0.5*\sc)};
	\draw[->,-stealth]  (2.55*\sc,-0.75*\sc) -- (2.4*\sc,-0.5*\sc);
	\node at (2.6*\sc,-0.9*\sc) {$\partial\Omega$};
	\node at (3.4*\sc,1.2*\sc) {$\mathbb{R}^n\setminus\Omega$};
	\end{tikzpicture}
	\caption{Domain of solving equations given in Section 2.} \label{Figure3}
\end{figure}

After introducing the microscopic models and deriving the macroscopic ones, we have insight into anomalous diffusions, especially L\'evy flights and tempered L\'evy flights.
In Section 2, all the derived equations are time dependent. From the process of derivation, one can see that the issue of initial condition can be easily/reasonably fixed, as classical ones, just specifying the value of $p({\bf X},0)$ in the domain $\Omega$.
For L\'evy processes, except Brownian motion, all others have discontinuous paths. As a result, the boundary $\partial \Omega$ itself (see Figure \ref{Figure3}) can not be hit by the majority of discontinuous sample trajectories. This implies that when solving the PDEs derived in Section 2, the generalized boundary conditions must be introduced, i.e., the information  of $p({\bf X},t)$ on the domain  ${\mathbb R}^n \backslash  \Omega$ must be properly accounted for.
In the following, we focus on Eqs. (\ref{FracLapl}), (\ref{HVequation2}), (\ref{temperedLaplacian}), (\ref{TemperedHVequation2}) to discuss the boundary issues.

\subsection{Generalized Dirichlet type boundary conditions}

The appropriate initial and boundary value problems for Eq. (\ref{FracLapl}) should be
\begin{equation}\label{IBFracLapl}\small
\left\{
\begin{aligned}
& \frac{\partial p({\bf X},t)}{\partial t}=\Delta^{\beta/2}p({\bf X},t)=\frac{-\beta \Gamma(\frac{n+\beta}{2})}{2^{1-\beta} \pi^{n/2}\Gamma(1-\beta/2)}\lim\limits_{\varepsilon \rightarrow 0^+} \int_{\mathscr C B_\varepsilon ({\bf X})} \frac{p({\bf X},t)-p({\bf Y},t)}{|{\bf X}- {\bf Y}|^{n+\beta}}d{\bf Y}\quad {\rm in}~~ \Omega, \\
& p({\bf X},0)|_{\Omega}=p_0({\bf X}), \\
& p({\bf X},t)|_{\mathbb{R}^n \backslash \Omega}=g({\bf X},t).
\end{aligned} \right.
\end{equation}
In Eq. (\ref{IBFracLapl}), the term
\begin{equation}\label{LaplTerm}
\begin{aligned}
&\lim\limits_{\varepsilon \rightarrow 0^+} \int_{\mathscr C B_\varepsilon ({\bf X})} \frac{p({\bf X},t)-p({\bf Y},t)}{|{\bf X}- {\bf Y}|^{n+\beta}}d{\bf Y} \\
& =\lim\limits_{\varepsilon \rightarrow 0^+} \int_{(\mathscr C B_\varepsilon ({\bf X})\cap\Omega)} \frac{p({\bf X},t)-p({\bf Y},t)}{|{\bf X}- {\bf Y}|^{n+\beta}}d{\bf Y}+\int_{\mathbb{R}^n\backslash\Omega} \frac{p({\bf X},t)-g({\bf Y},t)}{|{\bf X}- {\bf Y}|^{n+\beta}}d{\bf Y}
\\
& =\lim\limits_{\varepsilon \rightarrow 0^+} \int_{(\mathscr C B_\varepsilon ({\bf X})\cap\Omega)} \frac{p({\bf X},t)-p({\bf Y},t)}{|{\bf X}- {\bf Y}|^{n+\beta}}d{\bf Y}+p({\bf X},t)\int_{\mathbb{R}^n\backslash\Omega} {|{\bf X}- {\bf Y}|^{-n-\beta}}d{\bf Y}
\\
&~~~~+\int_{\mathbb{R}^n\backslash\Omega} \frac{-g({\bf Y},t)}{|{\bf X}- {\bf Y}|^{n+\beta}}d{\bf Y}.
\end{aligned}
\end{equation}
According to Eq. (\ref{LaplTerm}), $g({\bf X},t)$ should satisfy that there exist positive $M$ and $C$ such that when $|{\bf X}|>M$,
\begin{equation} \label{Requriments}
\frac{|g({\bf X},t)|}{|{\bf X}|^{\beta-\varepsilon}}<C~~ {\rm for~ positive~ small} ~ \varepsilon.
\end{equation}
In particular, when Eq. (\ref{Requriments}) holds, the function $\int_{\mathbb{R}^n\backslash\Omega} \frac{-g({\bf Y},t)}{|{\bf X}- {\bf Y}|^{n+\beta}}d{\bf Y}$ of ${\bf X}$ has any order of derivative if $g({\bf X},t)$ is integrable in any bounded domain. One of the most popular cases is $g({\bf X},t) \equiv 0$, which is the so-called absorbing boundary condition, implying that the particle is killed whenever it leaves the domain  $\Omega$. Another interesting case is for the steady state fraction diffusion equation
\begin{equation}\label{SteadyFracDiff}
\left\{
\begin{aligned}
& \Delta^{\beta/2}p({\bf X})=0~~ \text{in} ~ \Omega, \\
&p({\bf X})|_{\mathbb{R}^n \backslash \Omega}=g({\bf X}).
\end{aligned} \right.
\end{equation}
Given a domain $H \subset \mathbb{R}^n\backslash\Omega$, if taking $g({\bf X})=1$ for ${\bf X} \in H$ and $0$ for ${\bf X} \in (\mathbb{R}^n\backslash\Omega)\backslash H$, then the solution of (\ref{SteadyFracDiff}) means the probability that the particles undergoing L\'evy flights lands in $H$ after first escaping the domain $\Omega$ \cite{DengWW:2017}. If $g({\bf X}) \equiv 1$ in $\mathbb{R}^n\backslash\Omega$, then $p({\bf X})$ equals to $1$ in $\Omega$ because of the probability interpretation. This can also be analytically checked.

For the initial and boundary value problem Eq. (\ref{HVequation2}), it should be written as
\begin{equation} \label{IBHVequation2}
\left\{
\begin{aligned}
& \frac{\partial p({\bf x}_1,\cdots,{\bf x}_n,t)}{\partial t}=\frac{\partial^{\beta_1} p({\bf x}_1,\cdots,{\bf x}_n,t)}{\partial |{\bf x}_1|^{\beta_1} }+\frac{\partial^{\beta_2}  p({\bf x}_1,\cdots,{\bf x}_n,t)}{\partial |{\bf x}_2|^{\beta_2}}
\\
&\hspace*{3.1cm}+\cdots+\frac{\partial^{\beta_n}  p({\bf x}_1,\cdots,{\bf x}_n,t)}{\partial |{\bf x}_n|^{\beta_n}}
\quad  \mbox{in}\,\,\,\Omega, \\
&
p({\bf x}_1,\cdots,{\bf x}_n,0)|_{\Omega}=p_0({\bf x}_1,\cdots,{\bf x}_n), \\
&
p({\bf x}_1,\cdots,{\bf x}_n,t)|_{\mathbb{R}^n \backslash \Omega}=g({\bf x}_1,\cdots,{\bf x}_n,t).
\end{aligned}
\right.
\end{equation}
Similar to (\ref{LaplTerm}), in (\ref{IBHVequation2}) the term
\begin{equation}\label{HVLaplTerm}
\begin{aligned}
&\lim\limits_{\varepsilon \rightarrow 0^+} \int_{\mathscr C B_\varepsilon ({\bf x}_j)} \frac{p({\bf x}_1,\cdots,{\bf x}_j,\cdots,{\bf x}_n,t)-p({\bf x}_1,\cdots,{\bf y}_j,\cdots,{\bf x}_n,t)}{|{\bf x}_j- {\bf y}_j|^{1+\beta_j}}d{\bf y}_j
\\
& =\lim\limits_{\varepsilon \rightarrow 0^+} \int_{(\mathscr C B_\varepsilon ({\bf x}_j)\cap\Omega)} \frac{p({\bf x}_1,\cdots,{\bf x}_j,\cdots,{\bf x}_n,t)-p({\bf x}_1,\cdots,{\bf y}_j,\cdots,{\bf x}_n,t)}{|{\bf x}_j- {\bf y}_j|^{1+\beta_j}}d{\bf y}_j
\\
&
~~~+p({\bf x}_1,\cdots,{\bf x}_j,\cdots,{\bf x}_n,t)\int_{\mathbb{R}\backslash(\Omega \cap {\mathbb{R}_j})} {|{\bf x}_j- {\bf y}_j|^{-1-\beta_j}}d{\bf y}_j
\\
&
~~~+\int_{\mathbb{R}\backslash(\Omega \cap {\mathbb{R}_j})} \frac{-g({\bf x}_1,\cdots,{\bf y}_j,\cdots,{\bf x}_n,t)}{|{\bf x}_j- {\bf y}_j|^{1+\beta_j}}d{\bf y}_j.
\end{aligned}
\end{equation}

From Eq. (\ref{HVLaplTerm}), for $j=1,\cdots,n$, $g({\bf x}_1,\cdots,{\bf x}_j,\cdots,{\bf x}_n,t)$ should satisfies that there exist positive $M$ and $C$ such that when $|{\bf x}_j|>M$,
\begin{equation} \label{HVRequriments}
\frac{|g({\bf x}_1,\cdots,{\bf x}_j,\cdots,{\bf x}_n,t)|}{|{\bf x}_j|^{\beta_j-\varepsilon}}<C~~ {\rm for~ positive~ small} ~ \varepsilon.
\end{equation}
The discussions below Eq. (\ref{SteadyFracDiff}) still makes sense for Eq. (\ref{IBHVequation2}). If $g({\bf x}_1,\cdots,{\bf x}_j,\cdots,{\bf x}_n,t)$ satisfies Eq. (\ref{HVRequriments}), and it is integrable w.r.t. ${\bf x}_j$ in any bounded interval. Then $\int_{\mathbb{R}\backslash(\Omega \cap {\mathbb{R}_j})} \frac{-g({\bf x}_1,\cdots,{\bf y}_j,\cdots,{\bf x}_n,t)}{|{\bf x}_j- {\bf y}_j|^{1+\beta_j}}d{\bf y}_j$ has any order of partial derivative w.r.t. ${\bf x}_j$.

The initial and boundary value problem for Eq. (\ref{temperedLaplacian}) is
\begin{equation}\label{HVtemperedLaplacian}
\left\{
\begin{aligned}
	& \frac{\partial p({\bf X},t)}{\partial t}=(\Delta+\lambda)^{\beta/2}p({\bf X},t)
	\quad {\rm in}~~ \Omega, \\
	& p({\bf X},0)|_{\Omega}=p_0({\bf X}), \\
	& p({\bf X},t)|_{\mathbb{R}^n \backslash \Omega}=g({\bf X},t).
\end{aligned} \right.
\end{equation}
Like the discussions for Eq. (\ref{IBFracLapl}), $g({\bf X},t)$ should satisfies that there exist positive $M$ and $C$ such that when $|{\bf X}|>M$,
\begin{equation} \label{RequrimentsTempered}
\frac{|g({\bf X},t)|}{e^{(\lambda-\varepsilon){\bf |X|}}}<C~~ {\rm for~ positive~ small} ~ \varepsilon.
\end{equation}
 If Eq. (\ref{RequrimentsTempered}) holds and $g({\bf X},t)$ is integrable in any bounded domain, the function $\int_{\mathbb{R}^n\backslash\Omega} \frac{-g({\bf Y},t)}{e^{\lambda |{\bf X}- {\bf Y}|}|{\bf X}- {\bf Y}|^{n+\beta}}d{\bf Y}$ of ${\bf X}$ has any order of derivative.

Again, the corresponding tempered steady state fraction diffusion equation is
\begin{equation}\label{TemperedSteadyFracDiff}
\left\{
\begin{aligned}
& (\Delta+\lambda)^{\beta/2}p({\bf X})=0 \quad {\rm in ~~} \Omega, \\
&p({\bf X})|_{\mathbb{R}^n \backslash \Omega}=g({\bf X}).
\end{aligned} \right.
\end{equation}
For $H \subset \mathbb{R}^n\backslash\Omega$, if taking $g({\bf X})=1$ for ${\bf X} \in H$ and $0$ for ${\bf X} \in (\mathbb{R}^n\backslash\Omega)\backslash H$, then the solution of (\ref{TemperedSteadyFracDiff}) means the probability that the particles undergoing tempered L\'evy flights lands in $H$ after first escaping the domain $\Omega$. If $g({\bf X}) \equiv 1$ in $\mathbb{R}^n\backslash\Omega$, then $p({\bf X})$ equals to $1$ in $\Omega$.

The initial and boundary value problem  (\ref{TemperedHVequation2}) should be written as
\begin{equation} \label{temperedIBHVequation2}
\left\{
\begin{aligned}
& \frac{\partial p({\bf x}_1,\cdots,{\bf x}_n,t)}{\partial t}=\frac{\partial^{\beta_1,\lambda} p({\bf x}_1,\cdots,{\bf x}_n,t)}{\partial |{\bf x}_1|^{\beta_1,\lambda} }+\frac{\partial^{\beta_2,\lambda}  p({\bf x}_1,\cdots,{\bf x}_n,t)}{\partial |{\bf x}_2|^{\beta_2,\lambda}}
\\
&\hspace*{3.1cm}+\cdots+\frac{\partial^{\beta_n,\lambda}  p({\bf x}_1,\cdots,{\bf x}_n,t)}{\partial |{\bf x}_n|^{\beta_n,\lambda}}
\quad \mbox{in}\,\,\,\Omega,  \\
&
p({\bf x}_1,\cdots,{\bf x}_n,0)|_{\Omega}=p_0({\bf x}_1,\cdots,{\bf x}_n), \\
&
p({\bf x}_1,\cdots,{\bf x}_n,t)|_{\mathbb{R}^n \backslash \Omega}=g({\bf x}_1,\cdots,{\bf x}_n,t).
\end{aligned}
\right.
\end{equation}

For $j=1,\cdots,n$, $g({\bf x}_1,\cdots,{\bf x}_j,\cdots,{\bf x}_n,t)$ should satisfy that there exist positive $M$ and $C$ such that when $|{\bf x}_j|>M$,
\begin{equation} \label{TemperedHVRequriments}
\frac{|g({\bf x}_1,\cdots,{\bf x}_j,\cdots,{\bf x}_n,t)|}{e^{(\lambda-\varepsilon)|{\bf x}_j|}}<C~~ {\rm for~ positive~ small} ~ \varepsilon.
\end{equation}
If $g({\bf x}_1,\cdots,{\bf x}_j,\cdots,{\bf x}_n,t)$ is integrable w.r.t. ${\bf x}_j$ in any bounded interval and satisfies Eq. (\ref{TemperedHVRequriments}), then
$\int_{\mathbb{R}\backslash(\Omega \cap {\mathbb{R}_j})} \frac{-g({\bf x}_1,\cdots,{\bf y}_j,\cdots,{\bf x}_n,t)}{e^{\lambda|{\bf x}_j- {\bf y}_j}|{\bf x}_j- {\bf y}_j|^{1+\beta_j}}d{\bf y}_j$ has any order of partial derivative w.r.t. ${\bf x}_j$.

The ways of specifying the initial and boundary conditions for Eqs. (\ref{Feynman-Kac}) and (\ref{BackwardFeynman-Kac})
are the same as Eq. (\ref{IBFracLapl}). But for Eq. (\ref{Feynman-Kac}), the corresponding (\ref{Requriments}) should be changed as
\begin{equation} \label{Feynman-KacRequriments}
\frac{|U(\bf X)g({\bf X},t)|}{|{\bf X}|^{\beta-\varepsilon}}<C~~ {\rm for~ positive~ small} ~ \varepsilon.
\end{equation}
Similarly, the initial and boundary conditions of Eqs. (\ref{TemperedFeynman-Kac}) and (\ref{BackwardTemperedFeynman-Kac})
should be specified as the ones of Eq. (\ref{HVtemperedLaplacian}). But for Eq. (\ref{TemperedFeynman-Kac}), the corresponding (\ref{RequrimentsTempered}) needs to be changed as
\begin{equation} \label{Feynman-KacRequrimentsTempered}
\frac{|U({\bf X})g({\bf X},t)|}{e^{(\lambda-\varepsilon){\bf |X|}}}<C~~ {\rm for~ positive~ small} ~ \varepsilon.
\end{equation}
For the existence and uniqueness of the corresponding time-independent equations, one may refer to \cite{FelsingerKV:2015}.



\subsection{Generalized Neumann type boundary conditions}


Because of the inherent discontinuity of the trajectories of L{\'e}vy flights or tempered L{\'e}vy flights, the traditional Neumann type boundary conditions can not be simply extended to the fractional PDEs.  For the related discussions, see, e.g., \cite{BarlesGJ:2014,DipierroRV:2017}. Based on the models built in Sec. 2  and the law of mass conservation, we derive the reasonable ways of specifying the Neumann type boundary conditions, especially the reflecting ones. Let us first recall the derivation of classical diffusion equation. For normal diffusion (Brownian motion), microscopically the first moment of the distribution of waiting times and the second moment of the distribution of jump length are bounded, i.e., in Laplace and Fourier spaces, they are respectively like $1-c_1 u$ and $1-c_2 |{\bf k}|^2$; plugging them into Eq. (\ref{Montroll-WeissEq}) or Eq. (\ref{CharactFuncLevyFlight}) and performing integral transformations lead to the classical diffusion equation
\begin{equation}\label{ClassicalDE}
\frac{\partial p({\bf X},t)}{\partial t}=(c_2/c_1)\Delta p({\bf X},t).
\end{equation}
On the other hand, because of mass conservation, the continuity equation states that a change in density in any part of a system is due to inflow and outflow of particles into and out of that part of system, i.e., no particles are created or destroyed:
\begin{equation}\label{ContinuityEq}
\frac{\partial p({\bf X},t)}{\partial t}=-\nabla \cdot {\bf j},
\end{equation}
where ${\bf j}$ is the flux of diffusing particles. Combining (\ref{ClassicalDE}) with (\ref{ContinuityEq}), one may take
\begin{equation} \label{FicksLaw}
{\bf j}=-(c_2/c_1) \nabla p({\bf X},t),
\end{equation}
which is exactly Fick's law, a phenomenological postulation, saying that the flux goes from regions of high concentration to regions of low concentration with a magnitude proportional to the concentration gradient. In fact, for a long history, even up to now, most of the people are more familiar with the process: using the continuity equation (\ref{ContinuityEq})  and Fick's law (\ref{FicksLaw}) derives the diffusion equation (\ref{ClassicalDE}). The so-called reflecting boundary condition for (\ref{ClassicalDE}) is to let the flux ${\bf j}$ be zero along the boundary of considered domain.

 \begin{figure}
	\centering
	\begin{tikzpicture}[scale=1.5]
	\node at (1.5,0.5) {$\Omega$};
	\draw plot[smooth cycle] coordinates %
	{(0,0) (0.3,1) (1,1.5) (1.5, 1.6) (2,1.4) (2.5,0.95) (3,0) (2.0,-0.6) (1.5,-0.7) (0.5,-0.5)};
	\draw[->,-stealth] (-0.5,0.5) to [out=20,in=170] (3.0,1.4);
	\draw[->,-stealth] (3.2,1.0) to [out=240,in=30] (1.5,-1.0);
	\draw[->,-stealth] (0.6,-1.2) to [out=70,in=200] (2.2,0.3);
	\draw[->,-stealth] (1.3,0.6) to [out=220,in=50] (-0.3,-0.7);
	\end{tikzpicture}
	\caption{Sketch map of particles jumping into, or jumping out of, or passing through the domain: $\Omega$.}\label{Figure4}
\end{figure}
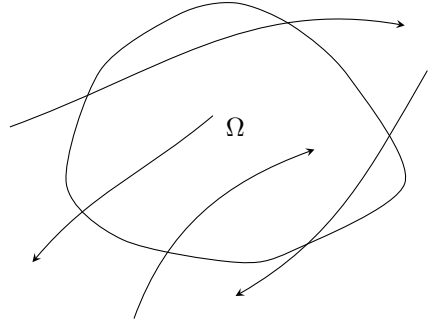

Here we want to stress that Eq. (\ref{ContinuityEq}) holds for any kind of diffusions, including the normal and anomalous ones. For  Eqs. (\ref{IBFracLapl},\ref{IBHVequation2},\ref{HVtemperedLaplacian},\ref{temperedIBHVequation2}) governing the PDF of L\'evy flights or tempered L\'evy flights, using the continuity equation  (\ref{ContinuityEq}), one can get the corresponding fluxes and the counterparts of Fick's law; may we call it fractional Fick's law.
Combining (\ref{IBFracLapl}) with (\ref{ContinuityEq}), one may let
\begin{equation}\label{LaplacianFlux}
{\bf j}_{\tiny\Delta}=\left\{-\frac{1}{2n} c_{n,\beta} \int_{-\infty}^{{\bf x}_i} \int_{\mathbb{R}^n} \frac{p({\bf X}+{\bf Y},t)+p({\bf X}-{\bf Y},t)-2 \cdot p({\bf X},t)}{|{\bf Y}|^{n+\beta}}d{\bf Y} d {\bf x}_i \right\}_{n \times 1}
\end{equation}
being the flux for the diffusion operator $\Delta^{\beta/2}$ with $\beta\in (0,2)$, or calling it fractional Fick's law corresponding to $\Delta^{\beta/2}$. From (\ref{IBHVequation2}) and (\ref{ContinuityEq}), one may choose
\begin{equation}\label{HVFlux}
{\bf j}_{\tiny {hv}}=\left\{-\frac{1}{2} c_{1,\beta_i} \int_{-\infty}^{{\bf x}_i} \int_{-\infty}^{+\infty} \frac{p({\bf X}+{\bf \widetilde{Y}}_i,t)+p({\bf X}-{\bf {\widetilde Y}}_i,t)-2 \cdot p({\bf X},t)}{|{\bf y}_i|^{1+\beta_i}}d{\bf y}_i d {\bf x}_i \right\}_{n \times 1},
\end{equation}
where ${\bf \widetilde{Y}}_i=\{{\bf x}_1,\dots, {\bf y}_i, \cdots,  {\bf x}_n \}^T$, being the flux (fractional Fick's law) corresponding to the horizontal and vertical type fractional operators. Similarly, we can also get the flux (fractional Fick's law) corresponding to the tempered fractional Laplacian and tempered horizontal and vertical type fractional operators, being  respectively taken as
\begin{equation}\label{TemperedLaplacianFlux}
{\bf j}_{\tiny\Delta,\lambda}=\left\{-\frac{1}{2n} c_{n,\beta, \lambda} \int_{-\infty}^{{\bf x}_i} \int_{\mathbb{R}^n} \frac{p({\bf X}+{\bf Y},t)+p({\bf X}-{\bf Y},t)-2 \cdot p({\bf X},t)}{e^{\lambda |{\bf Y}|}|{\bf Y}|^{n+\beta}}d{\bf Y} d {\bf x}_i \right\}_{n \times 1}
\end{equation}
and
\begin{equation}\label{TemperedHVFlux}
{\bf j}_{\tiny {hv}, \lambda}=\left\{-\frac{1}{2} c_{1,\beta_i,\lambda} \int_{-\infty}^{{\bf x}_i} \int_{-\infty}^{+\infty} \frac{p({\bf X}+{\bf \widetilde{Y}}_i,t)+p({\bf X}-{\bf {\widetilde Y}}_i,t)-2 \cdot p({\bf X},t)}{e^{\lambda |{\bf y}_i|}|{\bf y}_i|^{1+\beta_i}}d{\bf y}_i d {\bf x}_i \right\}_{n \times 1}
\end{equation}
with ${\bf \widetilde{Y}}_i=\{{\bf x}_1,\dots, {\bf y}_i, \cdots,  {\bf x}_n \}^T$.

Naturally, the Neumann type boundary conditions of (\ref{IBFracLapl},\ref{IBHVequation2},\ref{HVtemperedLaplacian},\ref{temperedIBHVequation2}) should be closely related to the values of the fluxes in the domain: $\mathbb{R}^n \backslash \Omega$;
if the fluxes are zero in it, then one gets the so-called reflecting boundary conditions of the equations. Microscopically the motion of particles undergoing L\'evy flights or tempered L\'evy flights are much different from the Brownian motion; very rare but extremely long jumps dominate the dynamics, making the trajectories of the particles discontinuous. As shown in Figure \ref{Figure4}, the particles may jump into, or jump out of, or even pass through the domain: $\Omega$. But the number of particles inside $\Omega$ is conservative, which can be easily verified by making the integration of (\ref{ContinuityEq}) in the domain $\Omega$, i.e.,
\begin{equation} \label{Conservation}
\frac{\partial }{\partial t} \int_{\Omega} p({\bf X},t) d {\bf X}=-\int_{\Omega} \nabla \cdot {\bf j} d {\bf X}=-\int_{\partial\Omega} {\bf j} \cdot {\bf n} ds=0,
\end{equation}
where ${\bf n}$ is the outward-pointing unit normal vector on the boundary. If ${\bf j}\,|_{\mathbb{R}^n \backslash \Omega}$=0, then for (\ref{IBFracLapl})
$ \Delta^{\frac{\beta}{2}} p({\bf X},t)=\nabla \cdot {\bf j} =0$  in $\mathbb{R}^n \backslash \Omega$. So, the Neumann type boundary conditions for   (\ref{IBFracLapl}), (\ref{IBHVequation2}), (\ref{HVtemperedLaplacian}), and (\ref{temperedIBHVequation2}) can be, heuristically, defined as
\begin{equation} \label{NeumannTypeBC1}
\Delta^{\frac{\beta}{2}} p({\bf X},t)=g({\bf X}) ~~ {\rm in} ~~  \mathbb{R}^n \backslash \Omega,
\end{equation}
\begin{equation} \label{NeumannTypeBC2}
\frac{\partial^{\beta_1} p({\bf x}_1,\cdots,{\bf x}_n,t)}{\partial |{\bf x}_1|^{\beta_1} }+\frac{\partial^{\beta_2}  p({\bf x}_1,\cdots,{\bf x}_n,t)}{\partial |{\bf x}_2|^{\beta_2}}+\cdots+\frac{\partial^{\beta_n}  p({\bf x}_1,\cdots,{\bf x}_n,t)}{\partial |{\bf x}_n|^{\beta_n}}=g({\bf X})
~~ {\rm in} ~~  \mathbb{R}^n \backslash \Omega,
\end{equation}
\begin{equation} \label{NeumannTypeBC3}
(\Delta+\lambda)^{\beta/2}p({\bf X},t)=g({\bf X})
~~ {\rm in} ~~  \mathbb{R}^n \backslash \Omega,
\end{equation}
and
\begin{equation} \label{NeumannTypeBC4}
\frac{\partial^{\beta_1,\lambda} p({\bf x}_1,\cdots,{\bf x}_n,t)}{\partial |{\bf x}_1|^{\beta_1,\lambda} }+\frac{\partial^{\beta_2,\lambda}  p({\bf x}_1,\cdots,{\bf x}_n,t)}{\partial |{\bf x}_2|^{\beta_2,\lambda}}+\cdots+\frac{\partial^{\beta_n,\lambda}  p({\bf x}_1,\cdots,{\bf x}_n,t)}{\partial |{\bf x}_n|^{\beta_n,\lambda}}=g({\bf X})
~~ {\rm in} ~~  \mathbb{R}^n \backslash \Omega,
\end{equation}
respectively. The corresponding reflecting boundary conditions are with $g({\bf X}) \equiv 0$.\medskip

{\bf Remark:}$\,\,$
The Neumann type boundary conditions \eqref{NeumannTypeBC1}-\eqref{NeumannTypeBC4} derived in this section are independent of the choice of the flux ${\bf j}$, provided that it satisfies
the condition \eqref{ContinuityEq}.

\section{Well-posedness and regularity of the fractional PDEs with generalized BCs}
Here, we show the well-posedesses of the models discussed in the above sections, taking the models with the operator $\Delta^{\frac{\beta}{2}}$ as examples; the other ones can be similarly proved.
For any real number $s\in{\mathbb R}$, we denote by $H^s({\mathbb R}^n)$ the conventional Sobolev space of functions  (see \cite{AdamsR:2003, McLean:2000}), equipped with the norm
$$
\|u\|_{H^s({\mathbb R}^n)}
:=\bigg(\int_{{\mathbb R}^n}(1+|{\bf k}|^{2s})|\widehat u({\bf k})|^2 d{\bf k}\bigg)^{\frac{1}{2}} ,
$$
The notation $H^s(\Omega)$ denotes the space of functions on $\Omega$ that admit extensions to $H^s({\mathbb R}^n)$, equipped with the quotient norm
$$
\|u\|_{H^s(\Omega)}
:=\inf_{\widetilde u}\|\widetilde u\|_{H^s({\mathbb R}^n)},
$$
where the infimum extends over all possible $\widetilde u\in H^s({\mathbb R}^n)$ such that $\widetilde u=u$ on $\Omega$ (in the sense of distributions).
The dual space of $H^s(\Omega)$ will be denoted by $H^s(\Omega)'$.
The following inequality will be used below:
\begin{align}\label{equiv-norm}
C^{-1}(\|\Delta^{\frac{\beta}{4}}u\|_{L^2({\mathbb R}^n)}
+ \|u\|_{L^2(\Omega)})
\le \|u\|_{H^{\frac{\beta}{2}}({\mathbb R}^n)}
\le C(\|\Delta^{\frac{\beta}{4}}u\|_{L^2({\mathbb R}^n)}
+ \|u\|_{L^2(\Omega)}) .
\end{align}

Let $H^s_0(\Omega)$ be the subspace of $H^s({\mathbb R}^n)$ consisting of functions which are zero  in ${\mathbb R}^n\backslash\Omega$.
It is isomorphic to the completion of $C^\infty_0(\Omega)$ in $H^s(\Omega)$.
The dual space of $H^s_0(\Omega)$ will be denoted by $H^{-s}(\Omega)$.

For any Banach space $B$, the space $L^2(0,T;B)$ consists of functions $u:(0,T)\rightarrow B$ such that
\begin{align}
\|u\|_{L^2(0,T;B)}
:=\bigg(\int_0^T\|u(\cdot,t)\|_{B}^2 d t\bigg)^{\frac{1}{2}} <\infty ,
\end{align}
and  $H^1(0,T;B)=\{u\in L^2(0,T;B):\partial_tu\in L^2(0,T;B)\}$; see \cite{Evans:2010}.

\subsection{Dirichlet problem}

For any given $g\in {\mathbb R} \cup (L^2(0,T;H^{\frac{\beta}{2}}({\mathbb R}^n))\cap H^1(0,T;H^{-\frac{\beta}{2}}({\mathbb R}^n))) \hookrightarrow C([0,T];L^2({\mathbb R}^n)$,
consider the time-dependent Dirichlet problem
\begin{align} \label{DrichletE1}
\left\{
\begin{aligned}
& \frac{\partial p}{\partial t}-\Delta^{\frac{\beta}{2}}p =f &&\mbox{in}\,\,\,\Omega , \\[3pt]
& p=g   &&\mbox{in}\,\,\,{\mathbb R}^n\backslash \Omega ,\\[5pt]
& p(\cdot,0)=p_0  &&{\mbox{in}\,\,\, \Omega ,}
\end{aligned} \right.
\end{align}
The weak formulation of \eqref{DrichletE1} is to find
$p=g+\phi$ such that
\begin{align}\label{PDE-reg}
\phi\in L^2(0,T;H^{\frac{\beta}{2}}_0(\Omega))\cap H^1(0,T;H^{-\frac{\beta}{2}}(\Omega))
\hookrightarrow C([0,T];L^2(\Omega))
\end{align}
and
\begin{align}\label{weak-phi-D}
 \int_{0}^{T}\int_\Omega \partial_t\phi \, q\, d {\bf X}d t
+  \int_{0}^{T}\int_{{\mathbb R}^n} \Delta^{\frac{\beta}{4}}\phi \,  \Delta^{\frac{\beta}{4}}q\,  d {\bf X} d t =
 \int_{0}^{T}\int_{\Omega} (f+\Delta^{\frac{\beta}{2}}g-\partial_tg)q\, d {\bf X} d t\\
\forall~ q\in L^2(0,T;H^{\frac{\beta}{2}}_0(\Omega)).\nonumber
\end{align}

It is easy to see that $a(\phi,q):=\int_{{\mathbb R}^n} \Delta^{\frac{\beta}{4}}\phi \,  \Delta^{\frac{\beta}{4}}q\,  d {\bf X} $ is a coercive bilinear form on $H^{\frac{\beta}{2}}_0(\Omega)\times H^{\frac{\beta}{2}}_0(\Omega)$ (cf. \cite[section 30.2]{Zeidler:1990})
and $\ell(q):=\int_{\Omega} (f+\Delta^{\frac{\beta}{2}}g-\partial_tg)q\, d {\bf X} $ is a continuous linear functional on $L^2(0,T;H^{\frac{\beta}{2}}_0(\Omega))$.
Such a problem as \eqref{weak-phi-D} has a unique weak solution
(cf. \cite[Theorem 30.A]{Zeidler:1990}).


The weak solution actually depends only on the values of $g$ in ${\mathbb R}^n\backslash\Omega$, independent of the values of $g$ in $\Omega$. To see this, suppose that  $g,\widetilde g\in {\mathbb R} \cup (L^2(0,T;H^{\frac{\beta}{2}}({\mathbb R}^n))\cap H^1(0,T;H^{-\frac{\beta}{2}}({\mathbb R}^n))) \hookrightarrow C([0,T];L^2({\mathbb R}^n))$ are two functions such that
$g=\widetilde g$ in ${\mathbb R}^n\backslash\Omega$, and
$p$ and $\widetilde p$ are the weak solutions of
\begin{align} \label{DrichletE}
\left\{
\begin{aligned}
& \frac{\partial p}{\partial t}-\Delta^{\frac{\beta}{2}}p =f &&\mbox{in}\,\,\,\Omega , \\[3pt]
& p=g   &&\mbox{in}\,\,\,{\mathbb R}^n\backslash \Omega ,\\[5pt]
& p(\cdot,0)=p_0  &&\mbox{in}\,\,\, \Omega,
\end{aligned} \right.
\qquad\mbox{and}\qquad
\left\{
\begin{aligned}
& \frac{\partial \widetilde p}{\partial t}-\Delta^{\frac{\beta}{2}}\widetilde p =f &&\mbox{in}\,\,\,\Omega , \\[3pt]
& \widetilde p=\widetilde  g   &&\mbox{in}\,\,\,{\mathbb R}^n\backslash \Omega ,\\[5pt]
& \widetilde p(\cdot,0)=p_0  &&\mbox{in}\,\,\, \Omega,
\end{aligned} \right.
\end{align}
respectively. Then the function $p-\widetilde p\in L^2(0,T;H^{\frac{\beta}{2}}_0(\Omega))\cap H^1(0,T;H^{-\frac{\beta}{2}}(\Omega)) $ satisfies
\begin{align}
& \int_{0}^{T}\int_\Omega \partial_t(p-\widetilde p)\, q\, d {\bf X}d t
+  \int_{0}^{T}\int_{{\mathbb R}^n} \Delta^{\frac{\beta}{4}}(p-\widetilde p) \,  \Delta^{\frac{\beta}{4}}q\,  d {\bf X} d t = 0
\quad
\forall\, q\in L^2(0,T;H^{\frac{\beta}{2}}_0(\Omega)) .
\end{align}
Substituting $q=p-\widetilde p$ into the equation above immediately yields
$p-\widetilde p=0$ a.e. in ${\mathbb R}^n\times(0,T)$.

\subsection{Neumann problem}
Consider the Neumann problem
\begin{align}\label{NeumannE}
\left\{
\begin{aligned}
&\frac{\partial p }{\partial t} - \Delta^{\frac{\beta}{2}}p =f
&&\mbox{in}\,\,\,\Omega , \\
&\Delta^{\frac{\beta}{2}}p =g
&&\mbox{in}\,\,\,{\mathbb R}^n\backslash\Omega,\\[5pt]
&p(\cdot,0)=p_0  &&\mbox{in}\,\,\,\Omega .
\end{aligned} \right.
\end{align}

\begin{definition}[Weak solutions]\label{weak-sol}
The weak formulation of \eqref{NeumannE} is to find $p\in L^2(0,T;H^{\frac{\beta}{2}}({\mathbb R}^n))\cap C([0,T];L^2(\Omega))$ such that
\begin{align}\label{PDE-reg1}
\partial_t p\in L^2(0,T;H^{\frac{\beta}{2}}(\Omega)') \text{ and } p(\cdot,0)=p_0,
\end{align}
satisfying the following equation:
\begin{equation}\label{PDE-ww}
  \begin{aligned}
& \int_{0}^{T}\int_\Omega \partial_tp({\bf X},t) q({\bf X},t)d {\bf X}d t
+  \int_{0}^{T}\int_{{\mathbb R}^n} \Delta^{\frac{\beta}{4}}p({\bf X},t)\Delta^{\frac{\beta}{4}}q({\bf X},t)d {\bf X} d t \\
&=
 \int_{0}^{T}\int_{\Omega} f({\bf X},t)q({\bf X},t)d {\bf X} d t
- \int_{0}^{T}\int_{{\mathbb R}^n\backslash\Omega}g({\bf X},t)q({\bf X},t)d {\bf X} d t
\\
&\forall\, q\in L^2(0,T;H^{\frac{\beta}{2}}({\mathbb R}^n)).
  \end{aligned}
\end{equation}
\end{definition}

\begin{theorem}[Existence and uniqueness of weak solutions]
If $p_0\in L^2(\Omega)$, $f\in L^2(0,T;H^{\frac{\beta}{2}}(\Omega)')$ and $g\in L^2(0,T;H^{\frac{\beta}{2}}({\mathbb R}^n\backslash\Omega)')$, then there exists a unique weak solution of \eqref{NeumannE} in the sense of Definition \ref{weak-sol}.
\end{theorem}

{\noindent\it Proof}$\,\,\,$
Let $t_k=k\tau$, $k=0,1,\dots,N$, be a partition of the time interval $[0,T]$, with step size $\tau=T/N$, and define
\begin{align}
&f_k({\bf X}):=\frac{1}{\tau}\int_{t_{k-1}}^{t_k} f({\bf X},t)d t ,\quad k=0,1,\dots,N,\\
&g_k({\bf X}):=\frac{1}{\tau}\int_{t_{k-1}}^{t_k} g({\bf X},t)d t ,\quad k=0,1,\dots,N .
\end{align}
Consider the time-discrete problem:
for a given $p_{k-1}\in L^2({\mathbb R}^n)$,
find $p_k\in H^{\frac{\beta}{2}}({\mathbb R}^n)$ such that the following equation holds:
\begin{align}\label{TD-PDE}
& \frac{1}{\tau} \int_\Omega p_k({\bf X})q({\bf X})d {\bf X}
+\int_{{\mathbb R}^n} \Delta^{\frac{\beta}{4}}p_k({\bf X})\Delta^{\frac{\beta}{4}}q({\bf X})d {\bf X} \nonumber \\
&=  \frac{1}{\tau} \int_\Omega p_{k-1}({\bf X})q({\bf X})d {\bf X}
+
\int_{\Omega} f_k({\bf X})q({\bf X})d {\bf X}
-\int_{{\mathbb R}^n\backslash\Omega}g_k({\bf X})q({\bf X})d {\bf X}
\quad\forall\, q\in H^{\frac{\beta}{2}}({\mathbb R}^n) .
\end{align}
In view of \eqref{equiv-norm}, the left-hand side of the equation above is a coercive bilinear form on $H^{\frac{\beta}{2}}({\mathbb R}^n)\times H^{\frac{\beta}{2}}({\mathbb R}^n)$, while the right-hand side is a continuous linear functional on
$H^{\frac{\beta}{2}}({\mathbb R}^n)$. Consequently, the Lax--Milgram Lemma implies that there exists a unique solution $p_k\in H^{\frac{\beta}{2}}({\mathbb R}^n)$ for \eqref{TD-PDE}.

Substituting $\displaystyle q= p_k$ into \eqref{TD-PDE} yields
\begin{align}
&\frac{\|p_k\|_{L^2(\Omega)}^2-\|p_{k-1}\|_{L^2(\Omega)}^2}{2\tau}
+ \|\Delta^{\frac{\beta}{4}}p_k\|_{L^2({\mathbb R}^n)}^2 \nonumber \\
&\le  \|f_k\|_{H^{\frac{\beta}{2}}(\Omega)'}\|p_k \|_{H^{\frac{\beta}{2}}(\Omega)} +\|g_k\|_{H^{\frac{\beta}{2}}({\mathbb R}^n\backslash\Omega)'}\|p_k \|_{H^{\frac{\beta}{2}}({\mathbb R}^n\backslash\Omega)}  \nonumber \\
&\le (\|f_k\|_{H^{\frac{\beta}{2}}(\Omega)'}+\|g_k\|_{H^{\frac{\beta}{2}}({\mathbb R}^n\backslash\Omega)'})
\|p_k \|_{H^{\frac{\beta}{2}}({\mathbb R}^n)} \nonumber \\
&\le (\|f_k\|_{H^{\frac{\beta}{2}}(\Omega)'}+\|g_k\|_{H^{\frac{\beta}{2}}({\mathbb R}^n\backslash\Omega)'})
( \|\Delta^{\frac{\beta}{4}}p_k\|_{L^2({\mathbb R}^n)}^2
+ \|p_k\|_{L^2(\Omega)}^2) .
\end{align}
Then, summing up the inequality above for $k=1,2,\dots,N$, we have
\begin{align}
&\max_{1\le k\le N} \|p_k\|_{L^2(\Omega)}^2
+ \tau \sum_{k=1}^N \|\Delta^{\frac{\beta}{4}}p_k\|_{L^2({\mathbb R}^n)}^2    \nonumber \\
&\le  \|p_0\|_{L^2(\Omega)}^2
+
C\tau \sum_{k=1}^N  (\|f_k\|_{H^{\frac{\beta}{2}}(\Omega)'}^2
+\|g_k\|_{H^{\frac{\beta}{2}}({\mathbb R}^n\backslash\Omega)'}^2
+\|p_k\|_{L^2(\Omega)}^2 )  .
\end{align}
By applying Gr\"onwall's inequality to the last estimate, there exists a positive constant $\tau_0$ such that when $\tau<\tau_0$ we have
\begin{align}\label{TD-Est00}
&\max_{1\le k\le N} \|p_k\|_{L^2(\Omega)}^2
+ \tau \sum_{k=1}^N \|p_k\|_{H^{\frac{\beta}{2}}({\mathbb R}^n)}^2    \nonumber \\
&\le C\|p_0\|_{L^2(\Omega)}^2
+
C\tau \sum_{k=1}^N  (\|f_k\|_{H^{\frac{\beta}{2}}(\Omega)'}^2
+\|g_k\|_{H^{\frac{\beta}{2}}({\mathbb R}^n\backslash\Omega)'}^2 )  .
\end{align}
Since any $q\in H^{\frac{\beta}{2}}(\Omega)$ can be extended to $q\in H^{\frac{\beta}{2}}({\mathbb R}^n)$ with
$\|q\|_{H^{\frac{\beta}{2}}({\mathbb R}^n)}\le 2\|q\|_{H^{\frac{\beta}{2}}(\Omega)} $,
choosing such a $q$ in \eqref{TD-PDE} yields
\begin{align*}
&\bigg|\int_\Omega \frac{p_k({\bf X})-p_{k-1}({\bf X})}{\tau} q({\bf X})d {\bf X}\bigg|\\
&=
\bigg|\int_{\Omega} f_k({\bf X})q({\bf X})d {\bf X}
-\int_{{\mathbb R}^n\backslash\Omega}g_k({\bf X})q({\bf X})d {\bf X}
-\int_{{\mathbb R}^n} \Delta^{\frac{\beta}{4}}p_k({\bf X})\Delta^{\frac{\beta}{4}}q({\bf X})d {\bf X} \bigg|\nonumber \\
&\le
C(\|f_k\|_{H^{\frac{\beta}{2}}(\Omega)'}
+\|g_k\|_{H^{\frac{\beta}{2}}({\mathbb R}^n\backslash\Omega)'}
+\|\Delta^{\frac{\beta}{4}}p_k\|_{L^2({\mathbb R}^n)} ) \|q\|_{H^{\frac{\beta}{2}}({\mathbb R}^n)} \\
&\le
C(\|f_k\|_{H^{\frac{\beta}{2}}(\Omega)'}
+\|g_k\|_{H^{\frac{\beta}{2}}({\mathbb R}^n\backslash\Omega)'}
+\|\Delta^{\frac{\beta}{4}}p_k\|_{L^2({\mathbb R}^n)} ) \|q\|_{H^{\frac{\beta}{2}}(\Omega)} .
\end{align*}
which implies (via duality)
\begin{align}
\bigg\|\frac{p_k-p_{k-1}}{\tau} \bigg\|_{H^{\frac{\beta}{2}}(\Omega)'}
&\le C (\|f_k\|_{H^{\frac{\beta}{2}}(\Omega)'}
+\|g_k\|_{H^{\frac{\beta}{2}}({\mathbb R}^n\backslash\Omega)'}
+\|\Delta^{\frac{\beta}{4}}p_k\|_{L^2({\mathbb R}^n)})  .
\end{align}
The last inequality and \eqref{TD-Est00} can be combined and written as
\begin{align}\label{TD-Est}
&\max_{1\le k\le N} \|p_k\|_{L^2(\Omega)}^2
+ \tau \sum_{k=1}^N
\bigg(\bigg\|\frac{p_k-p_{k-1}}{\tau} \bigg\|_{H^{\frac{\beta}{2}}(\Omega)'}^2
+\|p_k\|_{H^{\frac{\beta}{2}}({\mathbb R}^n)}^2 \bigg)   \nonumber \\
&\le C\|p_0\|_{L^2(\Omega)}^2
+
C\tau \sum_{k=1}^N  (\|f_k\|_{H^{\frac{\beta}{2}}(\Omega)'}^2
+\|g_k\|_{H^{\frac{\beta}{2}}({\mathbb R}^n\backslash\Omega)'}^2 )  .
\end{align}

If we define the piecewise constant functions
\begin{align}
&f^{(\tau)}({\bf X},t):=f_k({\bf X})=\frac{1}{\tau}\int_{t_{k-1}}^{t_k} f({\bf X},t)d t   && \mbox{for}\,\,\, t\in(t_{k-1},t_k],\,\,\, k=0,1,\dots,N,\\
&g^{(\tau)}({\bf X},t):=g_k({\bf X})=\frac{1}{\tau}\int_{t_{k-1}}^{t_k} g({\bf X},t)d t && \mbox{for}\,\,\, t\in(t_{k-1},t_k],\,\,\, k=0,1,\dots,N, \\
&p_+^{(\tau)}({\bf X},t):=p_k({\bf X})  &&\mbox{for}\,\,\, t\in(t_{k-1},t_k],\,\,\, k=0,1,\dots,N ,
\end{align}
and the piecewise linear function
\begin{align}
p^{(\tau)}({\bf X},t):=\frac{t_k-t}{\tau}p_{k-1}({\bf X})+\frac{t-t_{k-1}}{\tau}p_{k}({\bf X}) \quad \mbox{for}\,\,\, t\in[t_{k-1},t_k],\,\,\, k=0,1,\dots,N ,
\end{align}
then \eqref{TD-PDE} and \eqref{TD-Est} imply
\begin{equation}\label{TD-PDE-w}
  \begin{aligned}
& \int_{0}^{T}\int_\Omega \partial_tp^{(\tau)}({\bf X},t) q({\bf X},t)d {\bf X}d t
+  \int_{0}^{T}\int_{{\mathbb R}^n} \Delta^{\frac{\beta}{4}}p^{(\tau)}_+({\bf X},t)\Delta^{\frac{\beta}{4}}q({\bf X},t)d {\bf X} d t\nonumber \\
&=
 \int_{0}^{T}\int_{\Omega} f^{(\tau)}({\bf X},t)q({\bf X},t)d {\bf X} d t
- \int_{0}^{T}\int_{{\mathbb R}^n\backslash\Omega}g^{(\tau)}({\bf X},t)q({\bf X},t)d {\bf X} d t
\\
&\forall\, q\in L^2(0,T;H^{\frac{\beta}{2}}({\mathbb R}^n)),
  \end{aligned}
\end{equation}
and
\begin{equation}\label{TD-Est2}
  \begin{aligned}
&\|p^{(\tau)}\|_{C([0,T];L^2(\Omega))}+
\|\partial_tp^{(\tau)}\|_{L^2(0,T;H^{\frac{\beta}{2}}(\Omega)')}\\
&+\|p^{(\tau)}\|_{L^\infty(0,T;H^{\frac{\beta}{2}}({\mathbb R}^n))}
+\|p_+^{(\tau)}\|_{L^\infty(0,T;H^{\frac{\beta}{2}}({\mathbb R}^n))} \nonumber \\
&\le
C\left(\|f^{(\tau)}\|_{L^2(0,T;H^{\frac{\beta}{2}}(\Omega)')}
+\|g^{(\tau)}\|_{L^2(0,T;H^{\frac{\beta}{2}}({\mathbb R}^n\backslash\Omega)')}  \right) \nonumber \\
&\le
C\left(\|f\|_{L^2(0,T;H^{\frac{\beta}{2}}(\Omega)')}
+\|g\|_{L^2(0,T;H^{\frac{\beta}{2}}({\mathbb R}^n\backslash\Omega)')}  \right) ,
  \end{aligned}
\end{equation}
respectively,
where the constant $C$ is independent of the step size $\tau$.
The last inequality implies that $p^{(\tau)}$ is bounded in
$H^1(0,T;H^{\frac{\beta}{2}}(\Omega)')\cap L^2(0,T;H^{\frac{\beta}{2}}({\mathbb R}^n))
\hookrightarrow C([0,T];L^2(\Omega))$. Consequently, there exists
$p\in H^1(0,T;H^{\frac{\beta}{2}}(\Omega)')\cap L^2(0,T;H^{\frac{\beta}{2}}({\mathbb R}^n))\hookrightarrow C([0,T];L^2(\Omega))$ and a subsequence $\tau_j\rightarrow 0$ such that
\begin{align}
&\mbox{$p^{(\tau_j)}$ converges to $p$ weakly in $L^2(0,T;H^{\frac{\beta}{2}}({\mathbb R}^n)$}, \\
&\mbox{$p_+^{(\tau_j)}$ converges to $p$ weakly in $L^2(0,T;H^{\frac{\beta}{2}}({\mathbb R}^n)$}, \\
&\mbox{$\partial_tp^{(\tau_j)}$ converges to $\partial_tp$ weakly in $L^2(0,T;H^{\frac{\beta}{2}}(\Omega)')$},\\
&\mbox{$p^{(\tau_j)}$ converges to $p$ weakly in $C([0,T];H^{\frac{\beta}{2}}(\Omega)')$\,\,\, (see \cite[Appendix C]{Lions:1996}}) .
\end{align}
By taking $\tau=\tau_j\rightarrow 0$ in \eqref{TD-PDE-w}, we obtain \eqref{PDE-ww}.
This proves the existence of a weak solution $p$ satisfying \eqref{PDE-reg1}.

If there are two weak solutions $p$ and $\widetilde p$,
then their difference $\eta=p-\widetilde p$ satisfies the equation
\begin{align}
& \int_{0}^{T}\int_\Omega \partial_t(p-\widetilde p) q \, d {\bf X}d t
+  \int_{0}^{T}\int_{{\mathbb R}^n} \Delta^{\frac{\beta}{4}}(p-\widetilde p)\Delta^{\frac{\beta}{4}}q \, d {\bf X} d t =  0
\quad
\forall\, q\in L^2(0,T;H^{\frac{\beta}{2}}({\mathbb R}^n)) .
\end{align}
Substituting $q=p-\widetilde p$ into the equation yields
\begin{align}
& \| p(\cdot,t)-\widetilde p(\cdot,t)\|_{L^2(\Omega)}^2
+ \| \Delta^{\frac{\beta}{4}}(p-\widetilde p)\|_{L^2(0,T;L^2({\mathbb R}^n))}^2
=\| p(\cdot,0)-\widetilde p(\cdot,0)\|_{L^2(\Omega)}^2
=  0,
\end{align}
which implies $p=\widetilde p$ a.e. in ${\mathbb R}^n\times(0,T)$. The uniqueness is proved.
\qed\medskip

{\bf Remark:}$\,\,$
From the analysis of this section we see that,
although the initial data $p_0({\bf X})$ physically exists in the whole space ${\mathbb R}^n$,
one only needs to know its values in $\Omega$ to solve the PDEs (under both Dirichlet and Neumann boundary conditions).

\section{Conclusion}

In the past decades, fractional PDEs become popular as the effective models of characterizing L\'evy flights or tempered L\'evy flights. This paper is trying to answer the question: What are the physically meaningful and mathematically reasonable boundary constraints for the models? We physically introduce the process of the derivation of the fractional PDEs based on the microscopic models describing L\'evy flights or tempered L\'evy flights, and demonstrate that from a physical point of view when solving the fractional PDEs in a bounded domain $\Omega$, the informations of the models in  $\mathbb{R}^n\backslash\Omega$ should be involved. Inspired by the derivation process, we specify the Dirichlet type boundary constraint of the fractional PDEs as $p({\bf X},t)|_{\mathbb{R}^n\backslash \Omega}=g({\bf X},t)$ and Neumann type boundary constraints as, e.g., $(\Delta^{\beta/2}p({\bf X},t))|_{\mathbb{R}^n \backslash \Omega}=g({\bf X},t)$ for the fractional Laplacian operator.

The tempered fractional Laplacian operator  $(\Delta+\lambda)^{\beta/2}$ is physically introduced and mathematically defined. For the four specific fractional PDEs given in this paper, we prove their well-posedness with the specified Dirichlet or Neumann type boundary constraints. In fact, it can be easily checked that these fractional PDEs are not well-posed if their boundary constraints are (locally) given in the traditional way;  the potential reason is that locally dealing with the boundary contradicts with the principles that the L\'evy or tempered L\'evy flights follow.

\end{document}